\documentclass[reqno,a4paper,11pt]{amsart}
\usepackage{amsmath, amsthm, amscd, amsfonts, amssymb, graphicx, color}
\usepackage{amssymb}
\usepackage[english]{babel}
\usepackage[dvips, lmargin=4cm, rmargin=4cm, tmargin=4cm, bmargin=4cm]{geometry}
\usepackage[colorlinks=true,urlcolor=blue,linkcolor=blue, citecolor=red]{hyperref}
\usepackage{hyperref}
\usepackage[T1]{fontenc}
\usepackage{amssymb}
\usepackage{tikz}
\usetikzlibrary{backgrounds}
\usepackage[english]{babel}
\usepackage{enumitem}
\usepackage{dsfont}
\newtheorem{thm}{Theorem}[section]
\newtheorem{defini}{Definition}[section]
\newtheorem{rem}{Remark}[section]
\newtheorem{lem}{Lemma}[section]
\newtheorem{prop}{Proposition}[section]
\newtheorem{coro}{Corollary}[section]

\numberwithin{equation}{section}

\def \N {\mathbb{N} }

\def \R {\mathbb{R} }

\begin{document}
\title[Embedding and extension results]{Embedding and extension results in Fractional  Musielak-Sobolev spaces}
\author[E. Azroul, A. Benkirane, M. SHIMI  and  M. Srati]
{E. Azroul$^1$, A. Benkirane$^2$ M. SHIMI$^3$ and  M. Srati$^4$}
\address{E. Azroul, A. Benkirane, M Shimi and  M. Srati\newline
 Sidi Mohamed Ben Abdellah
 University,
 Faculty of Sciences Dhar El Mahraz, Laboratory of Mathematical Analysis and Applications, Fez, Morocco.}
\email{$^1$elhoussine.azroul@gmail.com}
\email{$^2$abd.benkirane@gmail.com}
\email{$^3$mohammed.shimi2@usmba.ac.ma}
\email{$^4$srati93@gmail.com}

\maketitle
\begin{abstract} 
In this paper,   we are concerned with some qualitative properties of the new fractional Musielak-Sobolev spaces $W^sL_{\varPhi_{x,y}}$ such that the generalized Poincar\'e type inequality and some continuous and compact embedding theorems of these spaces.
Moreover, we prove that any function in $W^sL_{\varPhi_{x,y}}(\Omega)$ may be extended to a
function in $W^sL_{\varPhi_{x,y}}(\R^N)$, with $\Omega \subset \R^N$ is a bounded domain  of class $C^{0,1}$. In addition, we establish a result  relates to the complemented subspace in $W^s{L_{\varPhi_{x,y}}}\left( \R^N\right)$. As
an application, using  the mountain pass theorem and some
variational methods, we investigate the existence of a nontrivial weak solution for a class of nonlocal fractional type problems with Dirichlet boundary data.
\end{abstract}
\section{Introduction and preliminaries results}\label{S1}

In the last decades, the fractional Sobolev spaces have been a classical topic in functional and harmonic analysis all along, see for instance \cite{NSL, ST} which treat the topic in detail. On the other hand, fractional spaces, and the corresponding nonlocal problems, are now experiencing impressive applications in different subjects, both for the pure mathematical research and in view of concrete applications, such as, phase transitions \cite{2}, crystal dislocation \cite{Ap8}, conservation laws \cite{C9},  mathematical finance \cite{26},  anomalous diffusion \cite{68}, stratified materials \cite{83}, etc.

As well known, when $s\in (0,1)$, $p\in [1, \infty)$ and $\Omega$ is an open subset of $\R^N$, with $C^{0,1}$-regularity, any function in fractional Sobolev space $W^{s,p}(\Omega)$ may be extended to a function in $W^{s,p}(\R^N)$ (see \cite{11}). Extension results are quite important in applications and for  improve certain embeddings theorems, in the integer classical case \cite{5,7,10,16,17}, as well as in the fractional case \cite{DD,11,9,Z}.  All these previous
results are held under certain crucial regularity assumptions on the domain $\Omega$.

The problem of characterizing the class of sets that are extension domains for $W^sL_{\varPhi_{x,y}}$ is open,  where $W^sL_{\varPhi_{x,y}}$ is the new fractional Musielak-Sobolev space introduced recently in \cite{benkirane}.  In general, an arbitrary open set $\Omega$ is not an extension domain for $W^sL_{\varPhi_{x,y}}$. When $s$ is an integer, we cite \cite{10} for a complete characterization in the special case when $s = 1$, $\varPhi_{x,y}(t)=t^2$ and $N=2$, and we refer the interested reader to the nice monograph of Leoni \cite{12} in which this problem is very well discussed (see, in particular, Chapter 11 and Chapter 12).\\

 The study of variational problems where the modular function  satisfies  nonpolynomial growth conditions instead of having the usual $p$-structure arouses much interest in the development of applications to electrorheological fluids as an important class of non-Newtonian fluids (sometimes referred to as smart fluids). The electro-rheological fluids are characterized by their ability to drastically change the mechanical properties under the influence of an external electromagnetic field. A mathematical model of electro-rheological fluids was proposed by Rajagopal and Ruzicka (we refer the reader to \cite{e1,e2} for more details). Another important application is related to image processing \cite{e3} where this kind of diffusion operator is used to underline the borders of the distorted image and to eliminate the noise.\\

Recently, great attention has been devoted to the study of a new class of fractional Sobolev spaces and related nonlocal problems, in particular, in the fractional Orlicz-Sobolev spaces $W^sL_\varPhi(\Omega)$ (see \cite{3,sr5,SRB,SR,bh2,sr_mo,sal1,sal2}) and in the fractional Sobolev spaces  with variable exponents $W^{s,p(x,y)}(\Omega)$ (see \cite{az,SH,SHG,SHN,SH2,SRH,Bh,ku}), in which the authors establish some basic properties of these modular spaces and  the associated nonlocal operators, they also obtained certain existence results for nonlocal problems involving this type of integro-differential operators. Furthermore,  in that context,  the authors in \cite{benkirane} introduced a new functional   framework  which can be seen as a natural generalization of the above mentioned functional spaces.\\

Our main purpose in this paper is to continue the study of the aforementioned new class of fractional Sobolev spaces and the related nonlocal operators, by introducing the new fractional Musielak-Sobolev spaces $W^sL_{\varPhi_{x,y}}(\Omega)$. In addition,  we establish certain basic properties of those spaces and we prove some continuous and compact embedding theorems of these spaces into Musielak spaces.
Moreover, we show that any function in fractional Musielak-Sobolev spaces $W^sL_{\varPhi_{x,y}}(\Omega)$ may be extended to a
function in $W^sL_{\varPhi_{x,y}}(\R^N)$, with $\Omega \subset \R^N$ is a bounded domain  of class $C^{0,1}$.  As
an application, we are concerned with the existence of  weak solutions of the following nonlocal problem
$$
 (\mathcal{P}_a) \hspace*{0.5cm} \left\{ 
    \begin{array}{clclc}
 
 (-\Delta)^{s_1}_{a_{(x,.)}} u + (-\Delta)^{s_2}_{a_{(x,.)}} u& = &  f(x,u)   & \text{ in }& \Omega, \\\\
    \hspace*{1cm} ~~~~~~u & = & 0 \hspace*{0.2cm} \hspace*{0.2cm} & \text{ in } & \R^N\setminus \Omega,
    \end{array}
    \right. 
\label{P} $$
  where $\Omega$ is an open bounded subset in $\R^N$, $N\geqslant 1$,   with Lipschitz boundary $\partial \Omega$, $0<s_2\leqslant s_1<1$, 
 $f: \Omega\times \R \longrightarrow \R$ is a Carath\'eodory function satisfies some suitable conditions. Moreover, for any $i=1,2$, $(-\Delta)^{s_i}_{a_{(x,.)}}$ is the nonlocal integro-differential operator of elliptic type defined as follows
   
              $$
               \begin{aligned}
               (-\Delta)^{s_i}_{a_{(x,.)}}u(x)=2\lim\limits_{\varepsilon\searrow 0} \int_{\R^N\setminus B_\varepsilon(x)} a_{(x,y)}\left( \dfrac{|u(x)-u(y)|}{|x-y|^{s_i} }\right)\dfrac{u(x)-u(y)}{|x-y|^{s_i}} \dfrac{dy}{|x-y|^{N+s_i} },
               \end{aligned}
                $$  
for all $x\in \R^N$, where $(x,y,t)\mapsto a_{(x,y)}(t):=a(x,y,t) : \overline{\Omega}\times\overline{\Omega}\times \R\longrightarrow \R$   is such that : $\varphi(.,.,.) : \overline{\Omega}\times\overline{\Omega}\times \R \longrightarrow \R$ defined by  
  $$
     \varphi_{x,y}(t):=\varphi(x,y,t)= \left\{ 
          \begin{array}{clclc}
        a(x,y,|t|)t   & \text{ for }& t\neq 0, \\\\
          0  & \text{ for } & t=0,
          \end{array}
          \right. 
       $$
is increasing homeomorphism from $\R$ onto itself. Let 
$$\varPhi_{x,y}(t):=\varPhi(x,y,t)=\int_{0}^{t}\varphi_{x,y}(\tau)d\tau~~\text{ for all } (x,y)\in \overline{\Omega}\times\overline{\Omega},~~\text{ and all } t\geqslant 0.$$  
Then, $\varPhi_{x,y}$ is a Musielak function (see \cite{mu}), that is
\begin{itemize}
\item  $\varPhi(x,y,.)$ is a $\varPhi$-function for every $(x,y)\in\overline{\Omega}\times\overline{\Omega}$, i.e.,   is continuous, nondecreasing function with $\varPhi(x,y,0)= 0$, $\varPhi(x,y,t)>0$ for $t>0$ and $\varPhi(x,y,t)\rightarrow \infty$ as $t\rightarrow \infty$.
\item For every $t\geqslant 0$, $\varPhi(.,.,t) : \overline{\Omega}\times\overline{\Omega} \longrightarrow \R$ is a measurable function.
\end{itemize}
We introduce the function $ \widehat{a}_x(t):=\widehat{a}(x,t)=a_{(x,x)}(t)$  for all $(x,t)\in \overline{\Omega}\times \R$. Then, we define the function $\widehat{\varphi}(.,.) : \overline{\Omega}\times \R \longrightarrow \R$  by  
  $$
     \widehat{\varphi}_{x}(t):=\widehat{\varphi}(x,t)= \left\{ 
          \begin{array}{clclc}
        \widehat{a}(x,|t|)t   & \text{ for }& t\neq 0, \\\\
          0  & \text{ for } & t=0,
          \end{array}
          \right. 
       $$
is increasing homeomorphism from $\R$ onto itself. If we set 
\begin{equation}\label{phi}
\widehat{\varPhi}_{x}(t):=\widehat{\varPhi}(x,t)=\int_{0}^{t}\widehat{\varphi}_{x}(\tau)d\tau ~~\text{ for all}~~ t\geqslant 0.
\end{equation}  
Then, $\widehat{\varPhi}_{x}$ is also a Musielak function.
\begin{rem}\text{}
	\begin{itemize}
\item	\textnormal {For all $(x,y)\in \overline{\Omega}\times\overline{\Omega}$,  $\varPhi_{x,y}$ and $\widehat{\varPhi}_{x}$ are two convex and increasing functions from $\R^+$ to $\R^+$.	
	\item Note that  $(-\Delta)^s_{a_{(x,.)}}$  is a nonlocal integro-differential operator of elliptic type which can be seen as a generalization of the fractional $p(x,.)$-Laplacian operator $(-\Delta_{p(x,.)})^{s}$ (when $a_{(x,y)}(t)=|t|^{p(x,y)-2}$, see for instance \cite{az,ku}), and of the fractional $p$-Laplacian operator $(-\Delta)_{p}^{s}$ in the constant exponent case ( when $p(x,y)=p=\text{constant}$, i.e. $a_{(x,y)}(t)=|t|^{p-2}$). 
	\item  The operator  $(-\Delta)^s_{a_{(x,.)}}$  reduces to the fractional $a$-Laplacian if $a_{(x,y)}(t)=a(t)$, i.e. the function $a$ is independent of variables $x,y$ (see for example \cite{3, sal1}). On the other hand, we remark that  is the fractional version of the well-known $a$-Laplacian operator $div\big(a(x,|\nabla u(x)|)\nabla u(x)\big)$ which is associated with the Musielak-Sobolev spaces (see \cite{ra}).\\
\item In contrast to the classical $p(x)$-Laplacian and $a$-Laplacian, which are  local operators, the integro-differential operator $(-\Delta)^s_{a_{(x,.)}}$  is a paradigm of the vast family of nonlocal nonlinear operators, and this has immediate consequences in the formulation of basic questions such as the Dirichlet problem. For this, the Dirichlet datum is given in $\mathbb{R}^{N}\setminus\Omega$ (which is different from the classical case of the $p(x)$-Laplacian and $a$-Laplacian) and not simply on $\partial\Omega$. Moreover, the value of $(-\Delta)^s_{a_{(x,.)}}u(x)$ at any point $x\in \Omega$ depends not only on the values of $u$ on $\Omega$, but actually on the entire space $\mathbb{R}^{N}$, which
implies that the first equation in \hyperref[P]{$(\mathcal{P}_{a})$} is no longer a pointwise equation, it  is no longer a pointwise identity. Hence, it is often called nonlocal problem. This causes
some mathematical difficulties which make the study of such a problem particularly interesting.}
\end{itemize}
\end{rem}
 This paper is organized as follows, In rest of this Section, we recall some preliminaries on the new fractional Musielak-Sobolev spaces which will be used frequently. In Section \ref{S2}, we establish some qualitative properties of these new spaces as the generalized Poincar\'e type inequality and some continuous and compact embedding theorems.
 Moreover, we prove that any function in fractional Musielak-Sobolev spaces $W^sL_{\varPhi_{x,y}}(\Omega)$ may be extended to a
 function in $W^sL_{\varPhi_{x,y}}(\R^N)$, with $\Omega \subset \R^N$ is a bounded domain  of class $C^{0,1}$. Furthermore,  we establish a result  relates to the complemented subspace in $W^s{L_{\varPhi_{x,y}}}\left( \R^N\right)$.  In Section \ref{S3}, as an application, we obtain the existence of a nontrivial weak solution  for problem \hyperref[P]{$(\mathcal{P}_{a})$},   by using a version of mountain pass theorem.\\ 

Next, For the function $\widehat{\varPhi}_x$ given in (\ref{phi}), we introduce the Musielak class as follows
 $$K_{\widehat{\varPhi}_x} (\Omega)=\left\lbrace u : \Omega \longrightarrow \R \text{ mesurable }: \int_\Omega\widehat{\varPhi}_x(|u(x)|)dx < \infty  \right\rbrace, $$
 and the Musielak space
  $$L_{\widehat{\varPhi}_x} (\Omega)=\left\lbrace u : \Omega \longrightarrow \R \text{ mesurable }: \int_\Omega\widehat{\varPhi}_x(\lambda |u(x)|)dx < \infty \text{ for some } \lambda>0 \right\rbrace. $$
The space $L_{\widehat{\varPhi}_x} (\Omega)$ is a Banach space endowed with the Luxemburg norm 
$$||u||_{\widehat{\varPhi}_x}=\inf\left\lbrace \lambda>0 \text{ : }\int_\Omega\widehat{\varPhi}_x\left( \dfrac{|u(x)|}{\lambda}\right) dx\leqslant 1\right\rbrace. $$
 The conjugate function of $\varPhi_{x,y}$ is defined by $\overline{\varPhi}_{x,y}(t)=\int_{0}^{t}\overline{\varphi}_{x,y}(\tau)d\tau~~\text{ for all } (x,y)\in\overline{\Omega}\times\overline{\Omega}~~ \text{ and all } t\geqslant 0$, where $\overline{\varphi}_{x,y} : \R\longrightarrow \R$ is given by $\overline{\varphi}_{x,y}(t):=\overline{\varphi}(x,y,t)=\sup\left\lbrace s \text{ : } \varphi(x,y,s)\leqslant t\right\rbrace.$ Furthermore, we have the following H\"older type inequality
  \begin{equation}
   \left| \int_{\Omega}uvdx\right| \leqslant 2||u||_{\widehat{\varPhi}_x}||v||_{\overline{\widehat{\varPhi}}_x}\hspace*{0.5cm} \text{ for all } u \in L_{\widehat{\varPhi}_x}(\Omega)  \text{ and } v\in L_{\overline{\widehat{\varPhi}}_x}(\Omega).
   \end{equation}
    Throughout this paper, we assume that
{\small$$\label{v1}
    1<\varphi^-:=\inf_{t\geqslant 0}\dfrac{t\varphi_{x,y}(t)}{\varPhi_{x,y}(t)}\leqslant \varphi^+:=\sup_{t\geqslant 0}\dfrac{t\varphi_{x,y}(t)}{\varPhi_{x,y}(t)}<+\infty\text{ for all } (x,y)\in\overline{\Omega}\times\overline{\Omega}.~~~~~~(\varPhi_1) $$}
    This relation implies  that
    \begin{equation}\label{A2}
        1<\widehat{\varphi}^-:=\inf_{t\geqslant 0}\dfrac{t\widehat{\varphi}_{x}(t)}{\widehat{\varPhi}_{x}(t)}\leqslant\widehat{\varphi}^+:=\sup_{t\geqslant 0}\dfrac{t\widehat{\varphi}_{x}(t)}{\widehat{\varPhi}_{x}(t)}<+\infty,\text{ for all } x\in\overline{\Omega}.\end{equation}
             It follows that  $\varPhi_{x,y}$ and $\widehat{\varPhi}_{x}$ satisfy the global $\Delta_2$-condition (see \cite{ra}), written $\varPhi_{x,y}\in \Delta_2$ and $\widehat{\varPhi}_{x}\in \Delta_2$, that is,
    \begin{equation}\label{r1}
    \varPhi_{x,y}(2t)\leqslant K_1\varPhi_{x,y}(t)~~ \text{ for all } (x,y)\in\overline{\Omega}\times\overline{\Omega},~~\text{ and  all } t\geqslant 0,
    \end{equation} and
    \begin{equation}\label{rr1}
        \widehat{\varPhi}_{x}(2t)\leqslant K_2\widehat{\varPhi}_{x}(t) ~~\text{ for any } x\in\overline{\Omega},~~\text{ and  all } t\geqslant 0,
        \end{equation}
 where $K_1$ and $K_2$ are two positive constants. The inequality $(\ref{rr1})$ implies that $L_{\widehat{\varPhi}_{x}}(\Omega)=K_{\widehat{\varPhi}_{x}}(\Omega)$ (see \cite{mu}).
 
  An important role in manipulating  the Musielak spaces is played by the modular of the space $L_{\widehat{\varPhi}_{x}}(\Omega)$. It is worth noticing that the relation between the norm and the modular shows an equivalence between the topology defined by the norm and that defined by the modular.
  \begin{prop}\label{pr1} \textnormal {(\cite[Proposition 2.1]{ra}).}
     Assume that the condition \hyperref[v1]{$(\varPhi_1)$} is satisfied. Then, for all $u \in L_{\widehat{\varPhi}_{x}}(\Omega)$, the following relations hold true 
     \begin{enumerate}[label=(\roman*)]
     	\item $||u||_{\widehat{\varPhi}_{x}}>1$ $\Rightarrow$ $||u||_{\widehat{\varPhi}_{x}}^{\widehat{\varphi}^-}\leqslant \int_{\Omega}\widehat{\varPhi}_{x}(|u|)dx\leqslant  ||u||_{\widehat{\varPhi}_{x}}^{\widehat{\varphi}^+},$
     	\item $||u||_{\widehat{\varPhi}_{x}}<1$ $\Rightarrow$ $||u||_{\widehat{\varPhi}_{x}}^{\widehat{\varphi}^+}\leqslant \int_{\Omega}\widehat{\varPhi}_{x}(|u|)dx\leqslant  ||u||_{\widehat{\varPhi}_{x}}^{\widehat{\varphi}^-}.$
     \end{enumerate}
    \end{prop}
    \begin{defini}
     Let $A_x(t)$, $B_x(t): \R^+\times \Omega\longrightarrow \R^+$ be two Musielak functions. 
            $A_x$ is stronger $($resp essentially stronger$)$ than $B_x$,  $A_x\succ B_x$ (resp $A_x\succ\succ B_x$) in symbols, if for almost every $x\in \overline{\Omega}$ 
          $$B(x,t)\leqslant A( x,a t), \text{   , } t\geqslant t_0\geqslant 0, $$
          for some $($resp for each$)$ $a>0$ and $t_0$ (depending on $a$).
    \end{defini}
            \begin{rem}$($\cite[Section 8.5]{1}$)$.
      $A_x\succ\succ B_x$  is equivalent to the condition \\
                     $$\lim_{t\rightarrow \infty}\left(\sup\limits_{x\in \overline{\Omega}}\dfrac{B( x,\lambda t)}{A(x,t)}\right) =0,$$
                     for all $\lambda>0$. 
                        \end{rem}
     
     \begin{thm}\label{2.1,}$($\cite{benkirane2}$)$.
               Let $\Omega$ be an open subset of $\R^N$ which has a finite volume. Let   $A_x(t)$, $B_x(t): ~\R^+\times \Omega\longrightarrow \R^+$ be two Musielak functions such that $B_x \prec\prec A_x$.Then any bounded subset $S$ of $L_{A_x}(\Omega)$ which is precompact  in $L^1(\Omega)$, is also precompact in $L_{B_x}(\Omega)$.
        \end{thm}

Now, due to the nonlocality of the operator $(-\Delta)^s_{a_{(x,.)}}$,  we  define the new fractional Musielak-Sobolev space as introduce in \cite{benkirane} as follows 
    \begingroup\makeatletter\def\f@size{9}\check@mathfonts$$ W^s{L_{\varPhi_{x,y}}}(\Omega)=\Bigg\{u\in L_{\widehat{\varPhi}_x}(\Omega) :  \int_{\Omega} \int_{\Omega} \varPhi_{x,y}\left( \dfrac{\lambda| u(x)- u(y)|}{|x-y|^s}\right) \dfrac{dxdy}{|x-y|^N}< \infty \text{ for some } \lambda >0 \Bigg\}.
$$\endgroup
This space can be equipped with the norm
\begin{equation}\label{r2}
||u||_{s,\varPhi_{x,y}}=||u||_{\widehat{\varPhi}_x}+[u]_{s,\varPhi_{x,y}},
\end{equation}
where $[.]_{s,\varPhi_{x,y}}$ is the Gagliardo seminorm defined by 
$$[u]_{s,\varPhi_{x,y}}=\inf \Bigg\{\lambda >0 :  \int_{\Omega} \int_{\Omega} \varPhi_{x,y}\left( \dfrac{|u(x)- u(y)|}{\lambda|x-y|^s}\right) \dfrac{dxdy}{|x-y|^N}\leqslant 1 \Bigg\}.
$$
\begin{rem}
\textnormal{ Fractinal Orlicz-Sobolev spaces (see e.g.  \cite{3,sr5,sal1}) and fractional Sobolev spaces with  variable exponent (see e.g.  \cite{az,Bh,ku}) are two distinct extensions
	of fractional Sobolev spaces $W^{s,p}(\Omega)$ (see  \cite{Bi,11}), and they are two special kinds of the new fractional Musielak-Sobolev spaces (see e.g. \cite{benkirane}). For the reader convenience, we recall the definition of the aforementioned spaces.}
\end{rem}
a$)-$ For the case: $\varPhi_{x,y}(t)=\varPhi(t)$, i.e. $\varPhi$ is independent of variables $x,y$, we say that $L_\varPhi$ and $W^sL_\varPhi$ are Orlicz spaces and fractional Orlicz-Sobolev  spaces respectively  such that
\begingroup\makeatletter\def\f@size{9.5}\check@mathfonts
$$
W^s{L_\varPhi}(\Omega)
      =\Bigg\{u\in L_\varPhi(\Omega) :  \int_{\Omega} \int_{\Omega} \varPhi\left( \dfrac{\lambda| u(x)- u(y)|}{|x-y|^s}\right)\dfrac{ dxdy}{|x-y|^N}< \infty \textnormal{ for some } \lambda >0 \Bigg\}.     
$$
       \endgroup
       
b$)-$ For the case: $\varPhi_{x,y}(t)=|t|^{p(x,y)}$ for all $(x,y)\in \overline{\Omega}\times\overline{\Omega}$,  where $p:\overline{\Omega}\times\overline{\Omega}\longrightarrow(1,+\infty)$ is a continuous bounded function such that
\begin{equation}\label{r4}
1<p^{-}=\underset{(x,y)\in \overline{\Omega}\times\overline{\Omega}}{\min}p(x,y)\leqslant p(x,y)\leqslant p^{+}=\underset{(x,y)\in \overline{\Omega}\times\overline{\Omega}}{\max}p(x,y)<+\infty,
\end{equation}
and
\begin{equation}\label{r3}
p ~\text{is symmetric, that is, }~~ p(x,y)=p(y,x)~~~ \text{for all }(x,y)\in\overline{\Omega}\times\overline{\Omega}.
\end{equation}
If denoted by
$\bar{p}(x)=p(x,x)$ for all $x\in \overline{\Omega}.$ Then,
 we replace $L_{\varPhi_x}$ by $L^{\overline{p}(x)}$, and $W^sL_{\varPhi_{x,y}}$ by $W^{s,p(x,y)}$ and we refer them as variable exponent Lebesgue spaces,  and fractional Sobolev spaces with variable exponent respectively, defined by
$$L^{\overline{p}(x)}(\Omega)=\bigg\{u:\Omega\longrightarrow \mathbb{R} ~~\text{measurable}: \int_{\Omega}|u(x)|^{\overline{p}(x)}dx<+\infty
\bigg\},$$ and 
 $$\hspace{-10cm}W=W^{s,p(x,y)}(\Omega)$$
\begingroup\makeatletter\def\f@size{10}\check@mathfonts $$\hspace*{0.6cm}=\bigg\{u\in L^{\bar{p}(x)}(\Omega): \int_{\Omega\times\Omega}\frac{|u(x)-u(y)|^{p(x,y)}}{\lambda^{p(x,y)}|x-y|^{sp(x,y)+N}}~dxdy <+\infty,~~ \text{for some}~~\lambda>0\bigg\}.$$\endgroup
with the norm
$$\|u\|_{W}=\|u\|_{L^{\bar{p}(x)}(\Omega)}+[u]_{W},$$
where $ [.]_{W}$ is a Gagliardo seminorm with variable exponent given by $$[u]_{W}=[u]_{s,p(x,y)}= \inf \bigg\{\lambda>0:\int_{\Omega\times\Omega}\frac{|u(x)-u(y)|^{p(x,y)}}{\lambda^{p(x,y)}|x-y|^{N+sp(x,y)}}~dxdy \leqslant1 \bigg\}.$$
When $p(x,y)= p= \text{constant} \in (1,+\infty)$, the space $W$ reduce to the fractional Sobolev space $W^{s,p}(\Omega)$. For  a comprehensive introduction to the study of these spaces and the related nonlocal problems, we refer the reader to \cite{Bi,11}.
 
\begin{thm}$($\cite{benkirane}$)$.
       Let $\Omega$ be an open subset of $\R^N$, and let $s\in (0,1)$. The space $W^sL_{\varPhi_{x,y}}(\Omega)$ is a Banach space with respect to the norm $(\ref{r2})$, and a  separable $($resp. reflexive$)$ space if and only if $\varPhi_{x,y} \in \Delta_2$ $($resp. $\varPhi_{x,y}\in \Delta_2 $ and $\overline{\varPhi}_{x,y}\in \Delta_2$$)$. Furthermore,
              if   $\varPhi_{x,y} \in \Delta_2$ and $\varPhi_{x,y}(\sqrt{t})$ is convex, then  the space $W^sL_{\varPhi_{x,y}}(\Omega)$ is an uniformly convex space.\end{thm}

           \begin{defini}$($\cite{benkirane}$)$.
           We say that $\varPhi_{x,y}$ satisfies the fractional boundedness condition, written $\varPhi_{x,y}\in \mathcal{B}_{f}$, if
           $$~~~~~~~~~~~~~~~~~~~~~~~~~~~~~~ \label{v2}         
           \sup\limits_{(x,y)\in \overline{\Omega}\times\overline{\Omega}}\varPhi_{x,y}(1)<\infty.~~~~~~~~~~~~~~~~~~~~~~~~~~(\varPhi_2)$$
           \end{defini}
           \begin{thm}  $($\cite{benkirane}$)$.    \label{TT}
                       Let $\Omega$ be an open subset of $\R^N$,  and  $0<s<1$. Assume that  $\varPhi_{x,y}\in \mathcal{B}_{f}$. 
                       Then,
                       $$C^2_0(\Omega)\subset W^sL_{\varPhi_{x,y}}(\Omega).$$
                  \end{thm}
       
    \begin{lem} $($\cite{benkirane}$)$.\label{ll}
    Assume that \hyperref[v1]{$(\varPhi_1)$} is satisfied. Then
    $$\overline{\varPhi}_{x,y}(\varphi_{x,y}(t))\leqslant \varphi^+\varPhi_{x,y}(t) ~~\text{ for all } (x,y) \in \overline{\Omega}\times\overline{\Omega}~~\text{ and all } t\geqslant 0.$$
    \end{lem} 
   
          \begin{lem}$($\cite{benkirane}$)$\label{2.2..}. Assume that \hyperref[v1]{$(\varPhi_1)$} is satisfied. Then the following inequalities hold true:
                   \begin{equation}\label{3.}
 \varPhi_{x,y}(\sigma t)\geqslant \sigma^{\varphi^-}\varPhi_{x,y}(t) ~~\text{ for all } t>0  \text{ and any  } \sigma>1,
                   \end{equation}
      \begin{equation}\label{3.2}
 \varPhi_{x,y}(\sigma t)\geqslant \sigma^{\varphi^+}\varPhi_{x,y}(t) ~~\text{ for all }  t>0  \text{ and any } \sigma\in (0,1),
                   \end{equation}
  \begin{equation}\label{r10}
                      \varPhi_{x,y}(\sigma t)\leqslant \sigma^{\varphi^+}\varPhi_{x,y}(t) ~~\text{ for all } t>0 \text{ and any } \sigma>1,
                      \end{equation} 
     \begin{equation}\label{r11}
                     \varPhi_{x,y}(t)\leqslant \sigma^{\varphi^-}\varPhi_{x,y}\left( \dfrac{t}{\sigma} \right) ~~\text{ for all } t>0 \text{ and any } \sigma \in(0,1).
                     \end{equation}                                  
                   \end{lem}  
                     For any $u \in W^sL_{\varPhi_{x,y}}(\Omega)$, we define the modular function on  $W^sL_{\varPhi_{x,y}}(\Omega)$  as follows  
                     $$\varPsi(u)=\displaystyle\int_{\Omega} \int_{\Omega} \varPhi_{x,y}\left( \dfrac{ |u(x)- u(y)|}{|x-y|^s}\right) \dfrac{dxdy}{|x-y|^N}+\int_{\Omega}\widehat{\varPhi}_{x}\left( |u(x)|\right) dx.$$

        \begin{prop}$($\cite{benkirane}$)$.\label{mod}
         Assume that \hyperref[v1]{$(\varPhi_1)$} is satisfied. Then, for any $u \in W^sL_{\varPhi_{x,y}}(\Omega)$, the following relations hold true:
           \begin{equation}\label{mod1}
     ||u||_{s,\varPhi_{x,y}}>1\Longrightarrow      ||u||_{s,\varPhi_{x,y}}^{\varphi^-} \leqslant  \varPsi(u)\leqslant  ||u||_{s,\varPhi_{x,y}}^{\varphi^+},
           \end{equation}
           \begin{equation}\label{mod2}
                ||u||_{s,\varPhi_{x,y}}<1\Longrightarrow    ||u||_{s,\varPhi_{x,y}}^{\varphi^+} \leqslant  \varPsi(u)\leqslant  ||u||_{s,\varPhi_{x,y}}^{\varphi^-}. \end{equation}
           \end{prop}  
                                               
 \begin{prop}\label{1111} $($\cite{benkirane}$)$ Suppose that \hyperref[v1]{$(\varPhi_1)$} is satisfied. Then, for any $u \in W^sL_{\varPhi_{x,y}}(\Omega)$, the following assertions hold true:                     
                    \begin{equation}\label{32}
                  [u]_{s,\varPhi_{x,y}}>1 \Longrightarrow  [u]^{\varphi^-}_{s,\varPhi_{x,y}}\leqslant \phi(u) \leqslant [u]^{\varphi^+}_{s,\varPhi_{x,y}},
                    \end{equation} 
                     \begin{equation}\label{A33}
                     [u]_{s,\varPhi_{x,y}}<1 \Longrightarrow  [u]^{\varphi^+}_{s,\varPhi_{x,y}}\leqslant \phi(u) \leqslant [u]^{\varphi^-}_{s,\varPhi_{x,y}}, 
                       \end{equation} 
         where $\phi (u)=  \displaystyle\int_{\Omega} \int_{\Omega} \varPhi_{x,y}\left( \dfrac{| u(x)- u(y)|}{|x-y|^s}\right)\dfrac{dxdy}{|x-y|^N}$.       
                    \end{prop}
   We conclude this section by recalling a version of the mountain pass theorem  as given in \cite{mont1} which is the main tool to prove our existence result.
       \begin{thm}\label{2.2}
       Let $X$ be a real Banach space and $J \in C^1(X,\R)$ with $J(0)=0$. Suppose that the following conditions hold:
       
       $(G_1)$ \label{G1} There exists 
       $\rho>0 \text{  and } r>0 \text{ such that } J(u)\geqslant r \text{ for } ||u||=\rho$.
        
       $(G_2)$ \label{G2} There exists
       $e \in X \text{ with } ||e||>\rho \text{  such that } J(e)\leqslant 0$.\\
       Let
       $$c:=\inf_{\gamma \in \Gamma} \max_{t\in [0,1]}J (\gamma (t)) \text{ with } \Gamma=\left\lbrace \gamma \in C([0,1],X); \gamma(0)=0 , \gamma(1)=e \right\rbrace.$$
       Then there exists a sequence $\left\lbrace u_n\right\rbrace $ in $X$ such that 
       $$J(u_n)\rightarrow c \text{ \hspace{1cm} and \hspace{1cm}} J'(u_n)\rightarrow 0.$$
        \end{thm}     
                                      
    \section{Embeddings and extension results}\label{S2}
  In this section, we will establish some embeddings and extension results of the new fractional Musielak Sobolev spaces.                
  \subsection{Some embedding results}$\label{101}$
   This subsection is devoted to the embedding results of the new fractional  Musielak-Sobolev spaces  $W^{s}L_{\varPhi_{x,y}}(\Omega)$. To this end, we follow the same approach used to obtain the embedding results in the fractional Orlicz-Sobolev space $W^{s}L_\varPhi(\Omega)$  
  established in \cite{3}.
  
  First, it is worth noticing that, as in the classical case with $s_1, s_2$ being an integers,
  the space $W^{s_1}L_{\varPhi_{x,y}}(\Omega)$ is continuously embedded in $W^{s_2}L_{\varPhi_{x,y}}(\Omega)$ when $s_2\leqslant s_1$, as next
  result points out.
  \begin{prop}\label{pro21}
  Let $\Omega$ be an open subset of $\R^N$ and let $0<s_2\leqslant s_1<1$. Assume that \hyperref[v1]{$(\varPhi_1)$}-\hyperref[v2]{$(\varPhi_2)$} holds true.  Then, $W^{s_1}{L_{\varPhi_{x,y}}}(\Omega)$ continuously embedded in $W^{s_2}{L_{\varPhi_{x,y}}}(\Omega)$.
  \end{prop}
              \noindent \textbf{Proof}.
      Let $u\in W^{s_1}{L_{\varPhi_{x,y}}}(\Omega)$ and $\lambda>0$, If we define 
       $$\Omega_1=\left\lbrace (x,y) \subset \Omega\times\Omega~~\ ~~ \dfrac{|D^{s_1}u|}{\lambda}\leqslant 1\right\rbrace \text{ and } \Omega_2=\Omega\times\Omega \setminus \Omega_1,$$
        where  $D^{s_i}u=\dfrac{ u(x)- u(y)}{|x-y|^{s_i}}, i=1,2,$ and we set $d\mu= \dfrac{dxdy}{|x-y|^N}$. 
        Then
   
     {\small$$  
     \begin{aligned}
     \int_{\Omega} \int_{\Omega} \varPhi_{x,y}\left( \dfrac{|D^{s_2}u|}{\lambda}\right)d\mu&=\int_{\Omega} \int_{\Omega} \varPhi_{x,y}\left( \dfrac{|D^{s_1}u|}{\lambda}\dfrac{1}{|x-y|^{s_2-s_1}}\right)d\mu\\
     &\leqslant \int_{\Omega} \int_{\Omega} \varPhi_{x,y}\left( \dfrac{|D^{s_1}u|}{\lambda}\right)\dfrac{ dxdy}{|x-y|^{N+p(s_2-s_1)}}\\
     &= \left( \int\int_{\Omega_1}+\int\int_{\Omega_2}\right) \varPhi_{x,y}\left(\dfrac{ |D^{s_1}u|}{\lambda}\right)\dfrac{ dxdy}{|x-y|^{N+p(s_2-s_1)}}\\
     &= I_{\Omega_1}+I_{\Omega_2}
     \end{aligned}
     $$}
   where $p=1$ if $|x-y|\geqslant 1$ or $p=\varphi^+$ if $|x-y|<1$ for all $(x,y)\in \Omega\times \Omega$.  Notice that
     \begingroup\makeatletter\def\f@size{8}\check@mathfonts   $$ I_{\Omega_1}=  \int\int_{\Omega_1} \varPhi_{x,y}\left( \dfrac{|D^{s_1}u|}{\lambda}\right)\dfrac{ dxdy}{|x-y|^{N+p(s_2-s_1)}}\leqslant \sup\limits_{(x,y)\in \overline{\Omega}\times\overline{\Omega}}\varPhi_{x,y}(1)\int\int_{\Omega_1} \dfrac{ dxdy}{|x-y|^{N+p(s_2-s_1)}}:=c_1.$$\endgroup
     Since $N>N+p(s_2-s_1)$, it follows that the above integral  is finite. On the other hand, we have   $$
     \begin{aligned}
     I_{\Omega_2}&\leqslant \int\int_{\Omega_2} \varPhi_{x,y}\left( \dfrac{|D^{s_1}u|}{\lambda}\right)\dfrac{dxdy}{|x-y|^{N+p(s_2-s_1)}}\\
     &\leqslant d^{p(s_1-s_2)}  \int\int_{\Omega_2} \varPhi_{x,y}\left( \dfrac{|D^{s_1}u|}{\lambda}\right)\dfrac{dxdy}{|x-y|^{N}}\\
     & \leqslant d^{p(s_1-s_2)}  \int_{\Omega}\int_{\Omega} \varPhi_{x,y}\left( \dfrac{|D^{s_1}u|}{\lambda}\right)d\mu<\infty,
     \end{aligned}$$
     where $d=\sup\limits_{(x,y)\in \overline{\Omega}\times\overline{\Omega}}|x-y|$. 
     Hence, 
     $$  
     \begin{aligned}
     \int_{\Omega} \int_{\Omega} \varPhi_{x,y}\left( \dfrac{|D^{s_2}u|}{[u]_{s_1,\varPhi_{x,y}}}\right)d\mu&\leqslant c_1+d^{p(s_1-s_2)}.
     \end{aligned}
     $$
     This fact implies that 
     \begin{equation}\label{6e}
     [u]_{s_2,\varPhi_{x,y}}\leqslant \left( c_1+d^{p(s_1-s_2)}\right)  [u]_{s_1,\varPhi_{x,y}}.
     \end{equation}         
   Therefore,
      $$||u||_{s_2,\varPhi_{x,y}}\leqslant c||u||_{s_1,\varPhi_{x,y}},$$
      where $c=(1+c_1+d^{p(s_1-s_2)})$.
       \hspace*{7.5cm $\Box$}
       \\

  Given $s\in (0,1)$ and let $\widehat{\varPhi}_{x}$ as  defined in (\ref{phi}). We denote by $\widehat{\varPhi}_{x}^{-1}$ the inverse function of $\widehat{\varPhi}_{x}$ which satisfies the following conditions:
    \begin{equation}\label{15}
    \int_{0}^{1} \dfrac{\widehat{\varPhi}_{x}^{-1}(\tau)}{\tau^{\frac{N+s}{N}}}d\tau<\infty~~ \text{ for all } x\in \overline{\Omega},
    \end{equation}
    
    \begin{equation}\label{16n}
    \int_{1}^{\infty} \dfrac{\widehat{\varPhi}_{x}^{-1}(\tau)}{\tau^{\frac{N+s}{N}}}d\tau=\infty ~~\text{ for all }x\in \overline{\Omega}.
    \end{equation}
   Note that, if $\varphi_{x,y}(t)=|t|^{p(x,y)-1}$, then (\ref{15}) holds precisely when $sp(x,y)<N$ for all $(x,y)\in \overline{\Omega}\times \overline{\Omega}$.\\
    If (\ref{16n}) is satisfied, we define the inverse  Musielak conjugate function of $\widehat{\varPhi}_x$ as follows
    \begin{equation}\label{17}
    (\widehat{\varPhi}^*_{x,s})^{-1}(t)=\int_{0}^{t}\dfrac{\widehat{\varPhi}_{x}^{-1}(\tau)}{\tau^{\frac{N+s}{N}}}d\tau.
    \end{equation}
     \begin{thm}\label{3.4}
   Let $\Omega$  be a bounded open
    subset of  $\R^N$ with $C^{0,1}$-regularity 
      and bounded boundary. If $(\ref{15})$ and  $(\ref{16n})$  hold, then 
   \begin{equation}\label{18}
    W^s{L_{\varPhi_{x,y}}}(\Omega)\hookrightarrow L_ {\widehat{\varPhi}^*_{x,s}}(\Omega).
   \end{equation}
  \end{thm}
   \begin{thm}\label{th2}
            Let $\Omega$  be a bounded open
              subset of  $\R^N$ and  $C^{0,1}$-regularity 
                with bounded boundary. If $(\ref{15})$ and  $(\ref{16n})$  hold, then the embedding
             \begin{equation}\label{27}
              W^s{L_{\varPhi_{x,y}}}(\Omega)\hookrightarrow L_{B_x}(\Omega),
             \end{equation}
             is compact for all $B_x\prec\prec \widehat{\varPhi}^*_{x,s}$.
             \end{thm}
  
  The proof will be carried out in a several lemmas. The first one establishes
  an estimate for the Musielak conjugate function $ \widehat{\varPhi}^*_{x,s}$ defined in (\ref{17}).
  \begin{lem}$\label{ll3.2}$
    Let $s\in (0,1)$, we assume that  $(\ref{15})$- $(\ref{16n})$ holds true and let $ \widehat{\varPhi}^*_{x,s}$ be defined by $(\ref{17})$. Then for all $s'\in (0,s)$, the following
    conclusions may be drawn.
    \begin{enumerate}

    \item $[ \widehat{\varPhi}^*_{x,s}(t)]^{\frac{N-s'}{N}}$ is a Musielak function, in particular,  $ \widehat{\varPhi}^*_{x,s}$ is a Musielak function.
    
    \item For every $\varepsilon >0$, there exists a constant $ K_\varepsilon>0$ such that for every $t$,
      \end{enumerate}
    \begin{equation}\label{19}
    [ \widehat{\varPhi}^*_{x,s}(t)]^{\frac{N-s'}{N}}\leqslant \dfrac{1}{2\varepsilon} \widehat{\varPhi}^*_{x,s}(t)+\dfrac{K_\varepsilon}{\varepsilon}t.
    \end{equation}  
    \end{lem}
   \noindent \textbf{Proof.} The proof is similar to \cite[Lemma 3]{3}.         
   \begin{lem}$\label{3.3}$
     Let $\Omega$ be an open subset of $\R^N$, and $0<s<1$. Let $f$ satisfies a
     Lipschitz-condition on $\R$ and $f(0)=0$, then,
     \begin{enumerate}
      \item[(i)] For every $u\in W^{s,1}_{loc}(\Omega)$, if $g(x)=f(|u(x)|)$, then $g\in W^{s,1}_{loc}(\Omega)$.
      \item[(ii)] For every $u\in W^s{L_{\varPhi_{x,y}}}(\Omega)$, if $g(x)=f(|u(x)|)$, then $g\in  W^s{L_{\varPhi_{x,y}}}(\Omega)$.
         \end{enumerate}
      In particular, for every $u\in W^{s,1}(\Omega)$, if $g(x)=f(|u(x)|)$, then  $g\in W^{s,1}(\Omega)$.   
         
     \end{lem}
     \noindent \textbf{Proof} \label{A3.4.} To prove Lemma \ref{3.3}, we follow the same approach as in \cite[Lemma 4]{3}.
    \begin{lem} \label{l2.3} 
  Let  $\Omega$ be a bounded open subset of $\R^N$ and let $0<s'<s<1$.  Assume that the condition  \hyperref[v1]{$(\varPhi_1)$}-\hyperref[v2]{$(\varPhi_2)$} are satisfied, then  the space  $W^sL_{\varPhi_{x,y}}(\Omega)$ is continuously embedded in $W^{s',q}(\Omega)$ for all $q\in [1,\varphi^-]$.
  \end{lem} 
   \noindent \textbf{Proof}
   By $(\ref{A2})$  there exist $c>0$ such that
  \begin{equation}\label{si}
 |t|^{\widehat{\varphi}^-}\leqslant c \widehat{\varPhi}_x(t)~~ t>1~~\text{ for all } x\in \overline{\Omega}. \end{equation}
   Indeed, form $(\ref{A2})$ we have
   $$\widehat{\varphi}^-\leqslant \dfrac{t \widehat{\varphi}_x(t)}{\widehat{\varPhi}_x(t)}~~\forall t>1,$$
   so, 
   $$ \widehat{\varphi}^-\left[ \ln(t)\right] '\leqslant \left[ \ln(\widehat{\varPhi}_x(t))\right]'~~\forall t>1,$$
   this implies 
   $$|t|^{\varphi^-}\leqslant \dfrac{\widehat{\varPhi}_x(t)}{\widehat{\varPhi}_x(1)}\leqslant c\widehat{\varPhi}_x(t),~~\forall t>1$$
   where $c=\dfrac{1}{\inf\limits_{x\in \Omega}\widehat{\varPhi}_x(1)}$, note that by definition of $\widehat{\varPhi}_x$ and \hyperref[v2]{$(\varPhi_2)$}, $0<c<\infty$.

   Then,  for $u\in W^sL_{\varPhi_{x,y}}(\Omega)$, we have
  $$
  \begin{aligned}
  \int_{\Omega} |u(x)|^{\widehat{\varphi}^-}dx&\leqslant\int_{\Omega\cap \left\lbrace |u|\leqslant 1\right\rbrace} |u(x)|^{\widehat{\varphi}^-}dx+\int_{\Omega\cap \left\lbrace |u|> 1\right\rbrace} |u(x)|^{\widehat{\varphi}^-}dx\\
  &\leqslant |\Omega|+ c\int_\Omega \widehat{\varPhi}_x(|u(x)|)dx.
    \end{aligned}
   $$ 
   Hence, 
   \begin{equation}\label{lab1}
   ||u||_{\widehat{\varphi}^-}\leqslant c_1 ||u||_{\widehat{\varPhi}_x},\end{equation}
   where $c_1=|\Omega|+ c$. On the other hand, similar to Proposition $\ref{pro21}$ we have 
  \begin{equation}\label{lab2} [u]_{s',\varphi^-}\leqslant c_2 [u]_{s,\varPhi_{x,y}}.\end{equation}
   Then, combining $(\ref{lab1})$ with $(\ref{lab2})$, we obtain 
       
                   $$ ||u||_{s',\varphi^-}\leqslant c [u]_{s,\varPhi_{x,y}}.$$ 
                   This completes the proof.        \hspace*{8.2cm $\Box$}
   
              \noindent \textbf{Proof of Theorem \ref{3.4}}.
                Let $0<s'<s<1$,  $\sigma(t)=[\widehat{\varPhi}^*_{x,s}(t)]^{\frac{N-s'}{N}}$
                and $u\in W^s{L_{\varPhi_{x,y}}}(\Omega)$, we suppose for the moment that $u$ is bounded on $\Omega$ and not equal to zero in $L_{\widehat{\varPhi}_{x}}(\Omega)$, then  $\displaystyle\int_\Omega \widehat{\varPhi}^*_{x,s}\left( \dfrac{|u(x)|}{\lambda}\right) dx$  decreases continuously from infinity to zero as  $\lambda$ increases from zero to infinity and according, assumes the value unity for some positive value $k$ of $\lambda$, thus 
          \begin{equation}\label{22n}
           \int_\Omega \widehat{\varPhi}^*_{x,s}\left( \dfrac{|u(x)|}{k}\right) dx =1 \text{ , } k=||u||_{\widehat{\varPhi}^*_{x,s}}.
          \end{equation}
                Let $f(x)=\sigma\left( \dfrac{u(x)}{k}\right) $. Using Lemma \ref{A3.4.} $u\in W^{s',1}(\Omega)$, and $\sigma$ is Lipschitz function (see Lemma $\ref{ll3.2}$), so from Lemma \ref{3.3} we have $f\in W^{s',1}(\Omega)$, and  since $N>s'$, then by \cite[Theroem 4.58]{DD}, one has,
                $$W^{s',1}(\Omega) \hookrightarrow L^{\frac{N}{N-s'}}(\Omega).$$ 
                Hence 
                $$||f||_{L^{\frac{N}{N-s'}}}\leqslant k_1 \left( ||f||_{L^1}+[f]_{s',1}\right), $$
                and by (\ref{22n}), we get
                $$1= \left( \int_\Omega \widehat{\varPhi}^*_{x,s}\left( \dfrac{|u(x)|}{k}\right) dx \right) ^{\frac{N-s'}{N}}=||f||_{L^{\frac{N}{N-s'}}}.$$
                This implies that
               \begin{equation}\label{23}
                \begin{aligned}
                   1 &\leqslant  k_1 \left( ||f||_{L^1}+[f]_{s',1}\right)\\
                   &= k_1 \left( \int_{\Omega}\sigma\left( \dfrac{u(x)}{k}\right) dx+\int_{\Omega}\int_\Omega\dfrac{ |f(x)-f(y)|}{|x-y|^{N+s'}}dxdy\right) \\
        &=  k_1 \left( \int_{\Omega}\sigma\left( \dfrac{u(x)}{k}\right) dx+\int_{\Omega}\int_\Omega\dfrac{\left| \sigma\left( \dfrac{u(x)}{k}\right) -\sigma\left( \dfrac{u(y)}{k}\right) \right| }{|x-y|^{N+s'}}dxdy\right) \\
                  & = k_1 I_1 +k_1 I_2.
                     \end{aligned}   
               \end{equation}
              From (\ref{19}) and since $L_{\widehat{\varPhi}_x}(\Omega)\hookrightarrow L^1(\Omega)$, then  for $\varepsilon =k_1$, we have
          \begin{equation}\label{24}
          k_1I_1 \leqslant \dfrac{1}{2} \int_\Omega \widehat{\varPhi}^*_{x,s}\left( \dfrac{|u(x)|}{k}\right) dx+\dfrac{k_\varepsilon}{k}\int_\Omega |u(x)|dx \leqslant \dfrac{1}{2} +\dfrac{k'_\varepsilon}{k}||u||_{\widehat{\varPhi}_x},
          \end{equation}
            On the other hand, as $\sigma$ is a Lipschitz function, then there exists $C>0$ such that 
           $$ k_1I_2\leqslant \dfrac{C}{k}\int_{\Omega}\int_\Omega\dfrac{ |u(x)-u(y)|}{|x-y|^{N+s'}}dxdy. $$
       Next, by Lemma \ref{A3.4.}, we have
          \begin{equation}\label{25}
          \int_{\Omega}\int_\Omega\dfrac{ |u(x)-u(y)|}{|x-y|^{N+s'}}dxdy\leqslant C'[u]_{s,\varPhi_{x,y}},
          \end{equation}
          thus, 
            \begin{equation}\label{26}
            k_1I_2\leqslant \dfrac{C}{k} C'[u]_{s,\varPhi_{x,y}}.
            \end{equation}
            We pose $k_3=Ck_1 C'$. Combining (\ref{24})-(\ref{26}), we obtain 
            $$1\leqslant  \dfrac{1}{2} +\dfrac{k'_\varepsilon}{k}||u||_{\widehat{\varPhi}_x}+\dfrac{k_3}{k}[u]_{s,\varPhi_{x,y}},$$
            this implies that,
            $$ \dfrac{k}{2} \leqslant k'_\varepsilon ||u||_{\widehat{\varPhi}_x} + k_3[u]_{s,\varPhi_{x,y}}.$$
           Hence, we obtain, 
            $$k=||u||_{\widehat{\varPhi}^*_{x,s}}\leqslant k_4||u||_{s,\varPhi_{x,y}},$$
            where $k_4= \max \left\lbrace  2k'_{\varepsilon} , 2 k_3\right\rbrace $.\\
             Now, for $u\in W^s{L_{\varPhi_{x,y}}}(\Omega)$ arbitrary, we define  
             $$ u_n(x)=\left \{
             \begin{array}{lcl}
            u(x)\hspace*{1cm} \text{ if } |u(x)|\leqslant n,\\
             n \text{ sgn } u(x) \text{ if }|u(x)|>n.
             \end{array}
             \right .$$         
       Then,   $\left\lbrace u_n\right\rbrace $ is bounded
       and  by Lemma \ref{3.3}, $u_n\in W^s{L_{\varPhi_{x,y}}}(\Omega)$. Moreover $$||u_n||_{\widehat{\varPhi}^*_{x,s}}\leqslant k_4||u_n||_{s,\varPhi_{x,y}}\leqslant k_4||u||_{s,\varPhi_{x,y}}.$$
              Let $\lim\limits_{n \rightarrow \infty}||u_n||_{\widehat{\varPhi}^*_{x,s}}=k$, then $ k\leqslant k_4||u||_{s,\varPhi_{x,y}}$. Using the  Fatou's Lemma, we get,
             $$ \int_\Omega \widehat{\varPhi}^*_{x,s}\left( \dfrac{|u(x)|}{k}\right) dx \leqslant \lim\limits_{n\rightarrow \infty}\int_\Omega \widehat{\varPhi}_{x,s}\left( \dfrac{|u_n(x)|}{||u_n||_{\widehat{\varPhi}^*_{x,s}}}\right) dx<1, $$
             Consequently, $u\in L_{\widehat{\varPhi}^*_{x,s}}(\Omega)$ and $||u||_{\widehat{\varPhi}^*_{x,s}}\leqslant k\leqslant k_4||u||_{s,\varPhi_{x,y}}.$        
    \hspace*{3cm$\Box$ }  
               
                      \noindent \textbf{Proof of Theorem $\ref{th2}$}.  
          From Lemma \ref{l2.3}, we get
          $$ W^s{L_{\varPhi_{x,y}}}(\Omega)\hookrightarrow W^{s',1}(\Omega)\hookrightarrow L^1(\Omega).$$
         By \cite[Theroem 4.58]{DD}, the latter embedding being compact. Since a bounded subset $S$ of $W^s{L_{\varPhi_{x,y}}}(\Omega)$ is also a bounded subset of $L_{\widehat{\varPhi}^*_{x,s}}(\Omega)$ and precompact in $L^1(\Omega)$. Hence,  by Theorem  \ref{2.1,}  it is precompact in $L_{B_x}(\Omega)$.
    \hspace*{3cm$\Box$ }

  
               Now, combining  Lemma \ref{l2.3} and \cite[Theroem 4.58]{DD}, we obtain the following results.
             \begin{coro}\label{3.1.}
             Let $\Omega$ be a bounded open subset of $\R^N$ with $C^{0,1}$-regularity and bounded boundary. Let $0<s'<s<1$ and let $\varPhi_{x,y}$ be a Musielak function satisfies the condition   \hyperref[v1]{$(\varPhi_1)$},
             we define  
             $$
          \varphi_{s'}^*= \hspace*{0.1cm} \left\{ 
                           \begin{array}{clclc}
                        \frac{N\varphi^-}{N-s'\varphi^-}\hspace*{0.5cm} \text{ if } N>s'\varphi^-
                          \\\\
                            \infty\hspace*{0.5cm}  \text{ if }  N \leqslant s'\varphi^-.
                            \label{eq1}
                           \end{array}
                           \right. 
                     $$
             \begin{itemize}
                     \item[$\bullet$] If $s'\varphi^-<N$, then 
                    $W^s{L_{\varPhi_{x,y}}}(\Omega) \hookrightarrow L^{q}(\Omega),$ for all $q\in [1,\varphi^*_{s'}]$
                    and the embedding 
                 $W^s{L_{\varPhi_{x,y}}}(\Omega) \hookrightarrow L^{q}(\Omega)$
                 is compact for all $q\in [1 , \varphi_{s'}^*)$.\\      
                  \item[$\bullet$]  If $s'\varphi^-=N$,  then 
                            $W^s{L_{\varPhi_{x,y}}}(\Omega) \hookrightarrow L^{q}(\Omega),$ for all $q\in [1,\infty]$
                            and the embedding 
                         $W^s{L_{\varPhi_{x,y}}}(\Omega) \hookrightarrow L^{q}(\Omega)$
                         is compact for all $q\in [1 , \infty)$.\\                
                            \item[$\bullet$]   If $s'\varphi^->N$,  then the embedding 
      $W^s{L_{\varPhi_{x,y}}}(\Omega) \hookrightarrow L^{\infty}(\Omega)$
                                           is compact.
                                             \end{itemize}    \end{coro} 
     Next, we introduce a closed linear subspace of $W^sL_{\varPhi_{x,y}}(\Omega)$ as follows
      $$W^s_0L_{\varPhi_{x,y}}(\Omega)=\left\lbrace u\in W^sL_{\varPhi_{x,y}}(\R^N) \text{ : } u=0 \text{ a.e in } \R^N\setminus \Omega \right\rbrace.$$
      Then, we have the following generalized  Poincar\'e type inequality. 
  \begin{thm}\label{pc}
                          Let $\Omega$ be a bounded open subset of  $\R^N$ with $C^{0,1}$-regularity 
                                and bounded boundary,  let $s\in (0,1)$.
                           Then there exists a positive constant $\gamma$ such that    
                 $$ ||u||_{\widehat{\varPhi}_x}\leqslant \gamma [u]_{s,\varPhi_{x,y}} \text{ for all   }  u \in  W^s_0L_{\varPhi_{x,y}}(\Omega).$$
                         \end{thm}  
 \noindent \textbf{Proof}. Let $\widehat{\sigma}(t)=[\widehat{\varPhi}_{x}(t)]^{\frac{N-s'}{N}}$ where $s < s'$
                 and $u\in W_0^s{L_{\varPhi_{x,y}}}(\Omega)$, then  similarly to Lemma \ref{3.3},  for   $f(x)=\widehat{\sigma}\left( \dfrac{u(x)}{k}\right) $ with $k=||u||_{\widehat{\varPhi}_x}$, we have  $f\in W_0^{s',1}(\Omega)$, and  since $N>s'$, then by \cite[Theroem 4.58]{DD}, one has,
                 $$W_0^{s',1}(\Omega) \hookrightarrow L^{\frac{N}{N-s'}}(\Omega).$$ 
                 It follows that
                 $$||f||_{L^{\frac{N}{N-s'}}}\leqslant k_5 [f]_{s',1}, $$
                 and 
                 $$1= \left( \int_\Omega \widehat{\varPhi}_{x}\left( \dfrac{|u(x)|}{k}\right) dx \right) ^{\frac{N-s'}{N}}=||f||_{L^{\frac{N}{N-s'}}}.$$
                 this implies that,
                \begin{equation}\label{23}
                 \begin{aligned}
                    1 &\leqslant  k_5  [f]_{s',1}\\
                    &= k_5 \int_{\Omega}\int_\Omega\dfrac{ |f(x)-f(y)|}{|x-y|^{N+s'}}dxdy \\
         &=  k_5 \int_{\Omega}\int_\Omega\dfrac{|\sigma(\dfrac{u(x)}{k})-\sigma(\dfrac{u(y)}{k})|}{|x-y|^{N+s'}}dxdy \\
                   & \leqslant\dfrac{k_5c}{k}\int_{\Omega}\int_\Omega\dfrac{ |u(x)-u(y)|}{|x-y|^{N+s'}}dxdy.\\
                   &\leqslant\dfrac{k_5c}{k} [u]_{s,\varPhi_{x,y}}.
                      \end{aligned}.   
                \end{equation}
        Thus       
                 $$ ||u||_{\widehat{\varPhi}_x}\leqslant \gamma [u]_{s,\varPhi_{x,y}}$$           where $\gamma=k_5c$.      \hspace*{10cm$\Box$ }
   \begin{rem}\text{}
   	\begin{itemize}
   		\item[$a)-$] As a trivial consequence of Theorem $\ref{pc}$, for a bounded open subset $\Omega$ of  $\R^N,$  
   		there exists a positive constant $\lambda_1$ such that, 
   		\begin{equation}\label{lambda1}
   		\int_\Omega\widehat{\varPhi}_x\left( u(x)\right) dx\leqslant \lambda_1\int_\Omega\int_\Omega\varPhi_{x,y}\left( \dfrac{|u(x)-u(y)|}{|x-y|^s}\right) \dfrac{dxdy}{|x-y|^N},
   		\end{equation}        
   		for all $u\in W_0^s{L_{\varPhi_{x,y}}}(\Omega)$.
   		\item [$b)-$] From Theorem $\ref{pc}$, we deduce that $[.]_{s,\varPhi_{x,y}}$ is a norm on
   		$W^s_0L_{\varPhi_{x,y}}(\Omega)$ which is equivalent to the norm  $||u||_{s,\varPhi_{x,y}}$.
   	\end{itemize}
            
   \end{rem} 
                               
 \subsection{Extending a $W^s{L_{\varPhi_{x,y}}}(\Omega)$ function to the whole of $\mathbb{R}^N$}
    As usual, for any $k\in \N$ and $\alpha \in (0,1]$, we say that $\Omega$ is of class $C^{k,\alpha}$ if there exists $M>0$ such that for any $x\in \partial \Omega$ there exists a ball $B=B_r(x)$, $r>0$, and an isomorphism $T : Q\rightarrow B$ such that
\begin{equation}\label{11.}
 T\in C^{k,\alpha}(\overline{Q}), \hspace*{0.2cm}T^{-1}\in C^{k,\alpha}(\overline{B}),\hspace*{0.2cm} T(Q_+)=B\cap \Omega,\hspace*{0.2cm} T(Q_0)=B\cap \partial \Omega
\end{equation}
    $$ \text{and } ||T||_{C^{k,\alpha}(\overline{Q})}+||T^{-1}||_{C^{k,\alpha}(\overline{B})}\leqslant M,$$
    where $$Q:=\left\lbrace x=(x',x_N)\in \R^{N-1}\times \R : |x'|<1 \text{ and } |x_N|<1 \right\rbrace,$$ 
     $$Q_+:=\left\lbrace x=(x',x_N)\in \R^{N-1}\times \R : |x'|<1 \text{ and } 0<x_N<1 \right\rbrace,$$ 
     $$ \text{ and } Q_0:=\left\lbrace x\in Q : x_N=0\right\rbrace.$$
     
     Given an open bounded domain  $\Omega \in \R^N$, For any $s\in (0,1)$ and any Musielak function $\varPhi_{x,y}$, we say that  an open set $\Omega$ is an extension domains for $W^sL_{\varPhi_{x,y}}(\Omega)$ if there exists a positive constant  $C=C(N,s,\Omega)$ such that for every function $u\in W^sL_{\varPhi_{x,y}}(\Omega)$ there exists $\widetilde{u}\in W^sL_{\varPhi_{x,y}}(\R^N)$ with $\widetilde{u}(x)=u(x)$ for all $x$ in  $\Omega$.\\
     
     Hence, our aim in this subsection is to show that any open bounded set $\Omega$ of class $C^{0,1}$ with
     bounded boundary is an extension domain for $W^s{L_{\varPhi_{x,y}}}$. To this end, we follow the same approach as in \cite{extension,11}.
     
     In this subsection we assume that $\varPhi_{x,y}\in \mathcal{B}_{f}$ and 
     $$ \label{v3}~~~~~~~~~~~~~~~~~~~~~~~~~~~~~~~~
       \lim\limits_{t\rightarrow \infty}\dfrac{|t|^{\varphi^+}}{\widehat{\varPhi}_{x,s}^*(kt)}=0 ~~\forall k>0.~~~~~~~~~~~~~~~~~~~~~~~~~~(\varPhi_3)
       $$
    	
We start with some preliminary lemmas, in which we will construct the
extension to the whole of $\R^N$ of a function $u$ defined on $\Omega$ in two separated
cases: when the function $u $ is identically zero in a neighborhood of the
boundary $\partial \Omega$ and when $\Omega$ coincides with the half-space $\R^N_+$.
\begin{lem}\label{lem24}
Assume that \hyperref[v1]{$(\varPhi_1)$}-\hyperref[v3]{$(\varPhi_3)$} holds. Then for all $u\in W^s{L_{\varPhi_{x,y}}}(\Omega)$
$$\int_{\Omega}\varPhi_{x,y}(|u(x)|)dx<\infty~~\forall y\in\overline{\Omega}.$$
\end{lem}
 \noindent \textbf{Proof}. First, seminary to $(\ref{si})$, we have
 $$\varPhi_{x,y}(t)\leqslant |t|^{\varphi^{+}}~~\forall t>1.$$
 Then, 
 $$
 \begin{aligned}
 \int_{\Omega}\varPhi_{x,y}(|u(x)|)dx&\leqslant \int_\Omega\varPhi_{x,y}(1)dx+\int_\Omega |u(x)|^{\varphi^+}\\
 &\leqslant \sup\limits_{(x,y)\in \overline{\Omega}\times\overline{\Omega}}\varPhi_{x,y}(1) |\Omega|+\|u\|^{\varphi^+}_{\varphi^+}
  \end{aligned}
 $$
 On the other hand, by \hyperref[v3]{$(\varPhi_3)$} we can used Theorem $\ref{th2}$ and we have 
 $$\|u\|_{\varphi^+}\leqslant c \|u\|_{s,\varPhi_{x,y}}.$$
 So, 
  $$
  \begin{aligned}
  \int_{\Omega}\varPhi_{x,y}(|u(x)|)dx&\leqslant
   \sup\limits_{(x,y)\in \overline{\Omega}\times\overline{\Omega}}\varPhi_{x,y}(1) |\Omega|+\|u\|_{\varphi^+}^{\varphi^+}\\
   &\leqslant \sup\limits_{(x,y)\in \overline{\Omega}\times\overline{\Omega}}\varPhi_{x,y}(1) |\Omega|+c\|u\|^{\varphi^+}_{s,\varPhi_{x,y}}\\
   &<\infty.
     \end{aligned}
  $$
 \begin{lem}\label{4.1}
 Let $\Omega$ be an open bounded subset of  $\R^N$, $s\in(0,1)$, and $u\in W^s{L_{\varPhi_{x,y}}}(\Omega)$. If there exists a compact subset $K \subset  \Omega$ such that $u \equiv 0$ in $\Omega \setminus K$, then the extension function $\widetilde{u}$ defined as  
 
 $$       \widetilde{u}(x)= \left \{
        \begin{array}{clclc}
       u(x) \hspace{1cm}  & if & x\in \Omega,\\\\
          0 \hspace*{1cm} & if & x \in \R^N\setminus \Omega,
        \end{array}
        \right .  
$$
 belongs to $W^s{L_{\varPhi_{x,y}}}(\R^N)$.
 \end{lem} 
 
 \noindent \textbf{Proof}. Let $u\in W^s{L_{\varPhi_{x,y}}}(\Omega)$. So, clearly $\widetilde{u}\in L_{\widehat{\varPhi}_{x}}(\R^N)$. On the other hand, since $ \widetilde{u}$ in $\R^N\setminus \Omega$, then for some $\lambda>0$, we have
  \begingroup\makeatletter\def\f@size{10}\check@mathfonts
$$
\begin{aligned}
 \int_{\R^N}\int_{\R^N}  \varPhi_{x,y}\left( \dfrac{|\widetilde{u}(x)-\widetilde{u}(y)|}{\lambda|x-y|^s}\right)\dfrac{dxdy}{|x-y|^N}
 & =  \int_{\Omega}\int_{\Omega}\varPhi_{x,y}\left( \dfrac{|u(x)-u(y)|}{\lambda|x-y|^s}\right)\dfrac{dxdy}{|x-y|^N}\\
 &\hspace*{0.2cm}+2 \int_{\R^N\setminus \Omega}\int_{\Omega}\varPhi_{x,y}\left( \dfrac{|u(x)|}{\lambda|x-y|^s}\right)\dfrac{dxdy}{|x-y|^N},
 \end{aligned}
 $$\endgroup
 where the first term in the right hand-side is finite. Moreover,   for any $y\in \R^N\setminus K$ we have
 $$
 \begin{aligned} \int_{\R^N\setminus \Omega}\int_{\Omega} & \varPhi_{x,y}\left( \dfrac{|u(x)|}{\lambda|x-y|^s}\right)\dfrac{dxdy}{|x-y|^N}\\
 &\leqslant  \int_{\R^N\setminus \Omega}\int_{\Omega} \varPhi_{x,y}\left( \dfrac{|u(x)|}{\lambda}\right)\dfrac{dxdy}{|x-y|^{sp+N}}\\
 &\leqslant \int_{\R^N\setminus \Omega}\int_{\Omega}\varPhi_{x,y}\left( \dfrac{\mathds{1}_K|u(x)|}{\lambda}\right) \dfrac{dxdy}{dis(y , \partial K)^{sp+N}}\\
 & \leqslant \int_{\Omega}\varPhi_{x,y}\left( \dfrac{|u(x)|}{\lambda}\right) dx  \int_{\R^N\setminus \Omega} \dfrac{dy}{dis(y,\partial K)^{sp+N}}.
 \end{aligned}
 $$
where $p=1$ if $|x-y|\geqslant 1$ or $p=\varphi^+$ if $|x-y|<1$ for some $(x,y)\in \Omega\times \Omega$.    Since $dis(\partial\Omega,\partial K)>0$ and $N+sp>N$.  Then, by Lemma $\ref{lem24}$ 
$$\int_{\R^N}\int_{\R^N}  \varPhi_{x,y}\left( \dfrac{|\widetilde{u}(x)-\widetilde{u}(y)|}{\lambda|x-y|^s}\right)\dfrac{dxdy}{|x-y|^N}<\infty.$$
   \hspace*{12.8cm $\Box$}
  \begin{lem}\label{4.2}
  Let $s\in (0,1)$ and let $\Omega$ be an open subset of  $\R^N$, symmetric with respect to
  the coordinate $x_N$, and consider the sets $\Omega_+=\left\lbrace x\in \Omega :x_N>0\right\rbrace $ and $\Omega_-=\left\lbrace x\in \Omega :x_N\leqslant 0\right\rbrace $. Let $u$ be a function in $W^s{L_{\varPhi_{x,y}}}(\Omega_+)$, we define 
    $$       \widetilde{u}(x)= \left \{
           \begin{array}{clclc}
          u(x',x_N) \hspace{1cm}  & if & x_N\geqslant 0,\\\\
             u(x',-x_N) \hspace*{1cm} & if & x_N<0.
           \end{array}
           \right .  
   $$
 Then $\widetilde{u}$ belongs to $W^s{L_{\varPhi_{x,y}}}(\Omega).$ 
  \end{lem}
  \noindent \textbf{Proof}.
 We set $\widehat{x}=(x',-x_N)$, then by splitting the integrals,  for some $\lambda>0$, we get
  $$
  \begin{aligned}
  \int_{\Omega}\widehat{\varPhi}_{x}\left( \dfrac{|\widetilde{u}(x)|}{\lambda}\right) dx &=\int_{\Omega_+}\widehat{\varPhi}_{x}\left( \dfrac{|u(x',x_N)|}{\lambda}\right)dx+\int_{\Omega_-}\widehat{\varPhi}_{x}\left( \dfrac{|u(x',-x_N)|}{\lambda}\right) dx\\
  &=\int_{\Omega_+}\widehat{\varPhi}_{x}\left( \dfrac{|u(x)|}{\lambda}\right)dx+\int_{\Omega_+}\widehat{\varPhi}_{x}\left( \dfrac{|u(\widehat{x}',\widehat{x}_N)|}{\lambda}\right) dx\\
  &=2\int_{\Omega_+}\widehat{\varPhi}_{x}\left( \dfrac{|u(x)|}{\lambda}\right)dx\\
  &<\infty.
    \end{aligned}
  $$
  On the other hand, for all $\lambda>0$, we have 
   \begingroup\makeatletter\def\f@size{9}\check@mathfonts
  $$
  \begin{aligned}
  & \int_{\Omega}  \int_{\Omega}\varPhi_{x,y}\left( \dfrac{|\widetilde{u}(x)-\widetilde{u}(y)|}{\lambda|x-y|^s}\right)\dfrac{dxdy}{|x-y|^N} =      \int_{\R^N\setminus\Omega_+}\int_{\R^N\setminus \Omega_+}\varPhi_{x,y}\left( \dfrac{|u(x',-x_N)-u(y',-y_N)|}{\lambda|x-y|^s}\right)\dfrac{dxdy}{|x-y|^N} \\
   & +2 \int_{\Omega_+}\int_{\R^N\setminus \Omega_+}\varPhi_{x,y}\left( \dfrac{|u(x)-u(y',-y_N)|}{\lambda|x-y|^s}\right)\dfrac{dxdy}{|x-y|^N}+\int_{\Omega_+}\int_{\Omega_+}\varPhi_{x,y}\left( \dfrac{|u(x)-u(y)|}{\lambda|x-y|^s}\right)\dfrac{dxdy}{|x-y|^N}.
   \end{aligned}
   $$\endgroup
  By changing variable $\widehat{x}=(x',-x_N)$ and $\widehat{y}=(y',-y_N)$ we get
   $$
   \begin{aligned}
   \int_{\Omega}\int_{\Omega}\varPhi_{x,y}\left( \dfrac{|\widetilde{u}(x)-\widetilde{u}(y)|}{\lambda|x-y|^s}\right)\dfrac{dxdy}{|x-y|^N}&= 4 \int_{\Omega_+}  \int_{\Omega_+}\varPhi_{x,y}\left( \dfrac{|u(x)-u(y)|}{\lambda|x-y|^s}\right)\dfrac{dxdy}{|x-y|^N}\\
   &<\infty.
      \end{aligned}$$
   This concludes the proof.
    \hspace*{8.5cm $\Box$}
   \begin{lem}\label{4.3}
   Let $\Omega$ be an open subset of $\R^N$ and $s\in (0,1)$. Let  $u\in W^s{L_{\varPhi_{x,y}}}(\Omega)$ and $\psi\in C^{0,1}(\Omega)$, $0\leqslant \psi\leqslant 1$. Then $\psi u \in W^s{L_{\varPhi_{x,y}}}(\Omega)$.
   \end{lem}
     \noindent \textbf{Proof}. Let $u\in W^s{L_{\varPhi_{x,y}}}(\Omega)$. Then, it is clear that $\psi u \in L_{\widehat{\varPhi}_{x}}(\Omega)$. In addition, since $|\psi|\leqslant 1$, it follows that $||\psi u ||_{L_{\widehat{\varPhi}_x}(\Omega)}\leqslant ||u ||_{L_{\widehat{\varPhi}_x}(\Omega)}$ . On the other hand, adding and subtracting the factor $\psi(x)u(y)$, we have 
      \begingroup\makeatletter\def\f@size{8}\check@mathfonts 
     $$
     \begin{aligned}
     &\int_{\Omega}\int_{\Omega}\varPhi_{x,y}\left( \dfrac{|\psi(x)u(x)-\psi(y)u(y)|}{\lambda|x-y|^s}\right)\dfrac{dxdy}{|x-y|^N}\\
     &\leqslant\int_{\Omega}\int_{\Omega} \varPhi_{x,y}\left( \dfrac{2|\psi(x)u(x)-\psi(x)u(y)|}{2\lambda|x-y|^s}+\dfrac{2|\psi(x)u(y)-\psi(y)u(y)|}{2\lambda|x-y|^s}\right)\dfrac{dxdy}{|x-y|^N}\\
     &\leqslant \dfrac{2^{\varphi^+}}{2}\left( \int_{\Omega}\int_{\Omega} \varPhi_{x,y}\left( \dfrac{|\psi(x)u(x)-\psi(x)u(y)|}{\lambda|x-y|^s}\right)\dfrac{dxdy}{|x-y|^N} +\int_{\Omega}\int_{\Omega}\varPhi_{x,y}\left( \dfrac{|\psi(x)u(y)-\psi(y)u(y)|}{\lambda|x-y|^s}\right)\dfrac{dxdy}{|x-y|^N}\right) \\
    & \leqslant 2^{\varphi^+-1}\left(  \int_{\Omega}\int_{\Omega}\varPhi_{x,y}\left( \dfrac{|u(x)-u(y)|}{\lambda|x-y|^s}\right)\dfrac{dxdy}{|x-y|^N} +\int_{\Omega}\int_{\Omega}\varPhi_{x,y}\left( \dfrac{|u(y)(\psi(x)-\psi(y))|}{\lambda|x-y|^s}\right)\dfrac{dxdy}{|x-y|^N}\right).
     \end{aligned}
     $$\endgroup
 Since $\psi \in C^{0,1}(\Omega)$, we obtain
 $$
      \begin{aligned}
 \int_{\Omega}\int_{\Omega} & \varPhi_{x,y}\left( \dfrac{|u(y)(\psi(x)-\psi(y))|}{\lambda|x-y|^s}\right)\dfrac{dxdy}{|x-y|^N}\\
&=\int_{\Omega}\int_{\Omega\cap\left\lbrace |x-y|\geqslant 1\right\rbrace }\varPhi_{x,y}\left( \dfrac{|u(y)(\psi(x)-\psi(y))|}{\lambda|x-y|^s}\right)\dfrac{dxdy}{|x-y|^N}\\
&\hspace*{0.2cm}+\int_{\Omega}\int_{\Omega\cap\left\lbrace |x-y|\leqslant 1\right\rbrace }\varPhi_{x,y}\left( \dfrac{u(y)(\psi(x)-\psi(y))}{\lambda|x-y|^s}\right)\dfrac{dxdy}{|x-y|^N}  \\
&\leqslant  \int_{\Omega}\int_{\Omega\cap\left\lbrace |x-y|\geqslant 1\right\rbrace }\varPhi_{x,y}\left( \dfrac{|u(y)|}{\lambda|x-y|^s}\right)\dfrac{dxdy}{|x-y|^N}\\
&\hspace*{0.2cm}+\int_{\Omega}\int_{\Omega\cap\left\lbrace |x-y|\leqslant 1\right\rbrace }\varPhi_{x,y}\left( \dfrac{L |u(y)||x-y|}{\lambda|x-y|^s}\right)\dfrac{dxdy}{|x-y|^N}  \\
& \leqslant C\int_{\Omega} \varPhi_{x,y}\left(\dfrac{ |u(y)|}{\lambda}\right)dy,
          \end{aligned}
   $$ 
  where $L$ denotes the Lipschitz constant of $\psi$ and\\ $$C= \displaystyle\int_{\Omega\cap\left\lbrace |z|\geqslant 1\right\rbrace} \dfrac{1}{|z|^{s+N}}dz+  \int_{\Omega\cap\left\lbrace |z|\leqslant 1\right\rbrace} \dfrac{L^\alpha}{|z|^{(s-1)+N}}dz $$
  with $\alpha=1$ if $L<1$ and $\alpha=\varphi^+$ if $L>1$.  
  Note that the two above integrals are finite. In fact, the kernel $|z|^{-(s+N)}$ is summable when $|z|\geqslant 1$ since $N+s>N$. Moreover, as $(s-1)+N<N$, then the kernel  $|z|^{-((s-1)+N)}$ is summable when $|z|\leqslant 1$. Therefore  $\psi u \in W^s{L_{\varPhi_{x,y}}}(\Omega)$.
   \hspace*{12.7cm $\Box$}
   
  Now, we are in position to state and prove the main theorem of this subsection.
  \begin{thm}\label{4.1.}
  Let $\Omega$ be an open subset of $\R^N$ with $C^{0,1}$-regularity and bounded boundary. Then, $W^s{L_{\varPhi_{x,y}}}(\Omega)$ is continuously embedded in $W^s{L_{\varPhi_{x,y}}}(\R^N)$, namely for any $u\in W^s{L_{\varPhi_{x,y}}}(\Omega)$ there exists $\widetilde{u} \in W^s{L_{\varPhi_{x,y}}}(\R^N)$.
  \end{thm}
   \noindent \textbf{Proof}. Since $\partial \Omega$ is compact, we can find a finite number of balls $B_j$ such that $\partial \Omega \subset \cup_{j=1}^{k}B_j$ and so we can write $\R^N=\cup_{j=1}^{k}B_j\cup (\R^N\setminus \partial \Omega)$.
   
 If we consider this covering, there exists a partition of unity related
 to it, that is, there exist $k + 1$ smooth functions $\psi_0, \psi_1,..., \psi_k$ such that $supp(\psi_0)\subset \R^N$, $supp(\psi_j)\subset B_j$ for all $j\in \left\lbrace 0,...,k\right\rbrace $ and $\sum\limits_{j=0}^{k}\psi_j=1$.
 Clearly,
 $$u=\sum\limits_{j=0}^{k}\psi_ju.$$
 
 From Lemma \ref{4.3}, we know that $\psi_0u\in W^s{L_{\varPhi_{x,y}}}(\Omega)$. Furthermore, since
 $\psi_0u=0$ in a neighborhood of $\partial \Omega$, we can extend it to the whole of $\R^N$
 by setting
  $$ \left\{ \begin{array}{clclc} 
    \widetilde{\psi_0u}(x)=
            \psi_0u \hspace{1cm}  & if & x\in \Omega,\\\\
               0 \hspace*{1cm} & if & x \in \R^N\setminus \Omega.
                       \label{eq1}
             \end{array}         
           \right.     
       $$               
 and $\widetilde{\psi_0u}\in W^s{L_{\varPhi_{x,y}}}(\Omega)$. 
 For all  $j\in \left\lbrace 1,...,k\right\rbrace$, let us consider $u_{|B_j\cap \Omega}$ and set
 $$v_j(y):=u(T_j(y)) \text{ for any } y \in Q_+,$$
 where $T_j : Q\rightarrow B_j$ is the isomorphism of class $C^{0,1}(\Omega)$ defined in $(\ref{11.})$. Note that such a $T_j$ exists because $\Omega$ is an open subset of class $C^{0,1}$.
  
 Next, we claim that $v_j \in W^s{L_{\varPhi_{x,y}}}(Q_+)$. Indeed, using the standard changing
 variable formula by setting $x = T_j(\widehat{x})$, we obtain
\begin{equation}\label{9}
 \begin{aligned}\int_{Q_+}\int_{Q_+} & \varPhi_{x,y}\left( \dfrac{|v_j(\widehat{x})-v_j(\widehat{y})|}{\lambda|\widehat{x}-\widehat{y}|^s}\right)\dfrac{d\widehat{x}d\widehat{y}}{|\widehat{x}-\widehat{y}|^N}\\
  & = \int_{Q_+}\int_{Q_+}  \varPhi_{x,y}\left( \dfrac{|u(T_j(\widehat{x}))-u(T_j(\widehat{y}))|}{\lambda|\widehat{x}-\widehat{y}|^s}\right)\dfrac{d\widehat{x}d\widehat{y}}{|\widehat{x}-\widehat{y}|^N}\\
  & = \int_{B_j\cap\Omega}  \int_{B_j\cap\Omega}\varPhi_{x,y}\left( \dfrac{|u(x)-u(y)|}{\lambda|T^{-1}_j(x)-T^{-1}_j(y)|^s}\right)\det(T_j^{-1})\dfrac{dxdy}{|T^{-1}_j(x)-T^{-1}_j(y)|^N}\\
  & \leqslant \bar{C} \int_{B_j\cap\Omega}  \int_{B_j\cap\Omega}\varPhi_{x,y}\left( \dfrac{|u(x)-u(y)|}{\lambda|x-y|^s}\right)\dfrac{dxdy}{|x-y|^N}.
  \end{aligned}
\end{equation}
  Hence $v_j \in W^s{L_{\varPhi_{x,y}}}(Q_+)$. Moreover, using Lemma $\ref{4.2}$, we can extend $v_j$ to all $Q$, then the extension $\widetilde{v_j} \in W^s{L_{\varPhi_{x,y}}}(Q)$.
  
    We set
 $$w_j(x):=\widetilde{v_j}(T_j^{-1}(x)) \text{ for any } x\in B_j.$$
 By arguing as above, since $T_j$ is bi-Lipschitz, it follows that $w_j\in W^s{L_{\varPhi_{x,y}}}(B_j)$. Note that $w_j \equiv u $ (and consequently $\psi_jw_j\equiv \psi_ju$) on $B_j\cap\Omega$.
 By definition $\psi_j w_j$ has compact support in $B_j$ and therefore, as done for
$\psi_0u$, we can consider the extension $\widetilde{\psi_jw_j}$ to all $\R^N$ in such a way that
 $\widetilde{\psi_jw_j}\in W^s{L_{\varPhi_{x,y}}}(\R^N)$.
 
   Finally, we set
   $$ \widetilde{u}= \widetilde{\psi_0 u}+ \sum_{j=1}^{k}\widetilde{\psi_j \omega_j}$$
   which is the extension of $u$ defined on all $\R^N$. By construction, it is clear that
   $\widetilde{u}|_{\Omega}=u$.
   \hspace*{11.2cm $\Box$}  
  \subsection{Complemented subspaces in $W^s{L_{\varPhi_{x,y}}}(\mathbb{R}^N)$}
  Let $(M,\|.\|_{M})$ be a Banach space of measurable functions on $\mathbb{R}^N$ and $A$ is a measurable subset of $\mathbb{R}^N$ with positive Lebesgue measure. We define the trace space $M_{|A}$ as follows$$M_{|A}:=\left\lbrace u:A\rightarrow\mathbb{R} \text{ such that }  \exists U\in M \text{  with } U_{|A}=u \text{ a.e. }\right\rbrace. $$
  This space is endowed with the norm 
  $$\|u\|_{M_{|A}}=\inf\left\lbrace\|U\|_M: U\in M, U_{|A}=u  \text{ a.e.} \right\rbrace.$$
  We denote by  $\mathcal{T}U=U|_{A}$ the trace operator on $A$. In particular, the aforementioned construction can be applies to the new fractional Musielak-Sobolev space $W^s{L_{\varPhi_{x,y}}}(\mathbb{R}^N)$.\\ Note that the definitions of $L_{\widehat{\varPhi}_{x}}(\mathbb{R}^N)$ and  $W^s{L_{\varPhi_{x,y}}}(\mathbb{R}^N)$ is analogous to $L_{\widehat{\varPhi}_{x}}(\Omega)$ and  $W^s{L_{\varPhi_{x,y}}}(\Omega)$, one just replace occurrence of $\Omega$ by $\mathbb{R}^N$.\\
  
  Next, let us consider the trace
   operator
$$ \mathcal{T} : W^s{L_{\varPhi_{x,y}}}\left( \R^N\right) \longrightarrow W^s{L_{\varPhi_{x,y}}}(\Omega)$$
$$\hspace*{3cm}\tilde{u}\longrightarrow \mathcal{T}\tilde{u}=\tilde{u}|_{\Omega}=u.$$
\begin{rem}\label{r2.2}
	On other words, an open domain $\Omega\subset \R^N $ is a $W^s{L_{\varPhi_{x,y}}}$-extension domain if there exists a continuous linear extension operator 
$$ \mathcal{E} : W^s{L_{\varPhi_{x,y}}} (\Omega)\longrightarrow W^s{L_{\varPhi_{x,y}}}\left( \R^N\right)$$
$$\hspace*{1.3cm}u\longrightarrow \mathcal{E}u=\tilde{u}$$
\end{rem}
 \begin{defini}
 	A closed subspace
 	$Y$ of a Banach space $X$ is complemented if there is another closed subspace $Z$
 	of $X$ such that $X = Y \oplus Z$. That is, $Y \cap Z ={0}$ and every element $x\in X$
 	can be written as $x = y + z$, with $y\in Y$ and $z\in Z$.
 \end{defini} 
 In this subsection, our  result  relates to the complemented subspace in $W^s{L_{\varPhi_{x,y}}}\left( \R^N\right)$. More precisely, we have the following theorem,
\begin{thm}
Suppose that $\Omega\subset \R^N $ is an open bounded set of class $C^{0,1}$. Then
$$ W^s{L_{\varPhi_{x,y}}}(\R^N)=\ker\mathcal{T}\oplus \mathcal{E}(W^s{L_{\varPhi_{x,y}}}(\Omega)).$$

\end{thm}
 \noindent \textbf{Proof}. From the Theorem \ref{4.1.} and Remark \ref{r2.2}, we know that there exists a linear operator 
 $ \mathcal{E} : W^s{L_{\varPhi_{x,y}}}(\Omega)\longrightarrow W^s{L_{\varPhi_{x,y}}}(\R^N)$ such that
 \begin{enumerate}
 \item $\mathcal{E} u_{|\Omega}=u$,
 \item $\|\mathcal{E} u\|_{W^s{L_{\varPhi_{x,y}}}(\R^N)}\leqslant C||u||_{W^s{L_{\varPhi_{x,y}}}(\R^N)}.$
 \end{enumerate}
 On the one hand,  $\mathcal{E}(W^s{L_{\varPhi_{x,y}}}(\Omega))$ is a closed subspace of $W^s{L_{\varPhi_{x,y}}}(\R^N)$. On the other hand, for every element $u$ in $W^s{L_{\varPhi_{x,y}}}(\R^N)$ can be written as $u=u-\mathcal{E}(\mathcal{T}(u))+\mathcal{E}(\mathcal{T}(u))$, as $u-\mathcal{E}(\mathcal{T}(u)) \in \ker \mathcal{T}$ and $\varepsilon(\mathcal{T}(u)) \in \mathcal{E}(W^s{L_{\varPhi_{x,y}}}(\Omega))$. Moreover, $\ker\mathcal{T}\cap\mathcal{E}(W^s{L_{\varPhi_{x,y}}}(\Omega))={0}$, we deduce that
$$ W^s{L_{\varPhi_{x,y}}}(\R^N)=\ker\mathcal{T}\oplus \mathcal{E}(W^s{L_{\varPhi_{x,y}}}(\Omega)).$$
  \hspace*{12.7cm $\Box$}
  \section{Application to a nonlocal problem}\label{S3}
 In this section,  we investigate the existence of weak solution for  problem \hyperref[P]{$(\mathcal{P}_{a})$} in the new fractional Musielak-Sobolev space.

 We start by considering the function  $\alpha_x(t):=\alpha(x,t) : \overline{\Omega}\times \R \longrightarrow \R$ such that
 $g(.,.) : \overline{\Omega}\times \R \longrightarrow \R$ defined by  
   $$
      g_{x}(t):=g(x,t)= \left\{ 
           \begin{array}{clclc}
         \alpha(x,|t|)t   & \text{ for }& t\neq 0, \\\\
           0  & \text{ for } & t=0,
           \end{array}
           \right. 
        $$
 which is an increasing homeomorphism from $\R$ onto itself. We assume that $f: \Omega\times \R \longrightarrow \R$ is a Carath\'eodory function satisfies the following condition    $$ (f_1):\hspace*{3.8cm}\label{f1}  |f(x,t)|\leqslant c_0(1 + g_x(|t|)),  \hspace*{3.8cm} $$  
       where  $c_0$ is a nonnegative constant. Furthermore, if we set 
              \begin{equation}\label{10}
               G_x(t):=G(x,t)=\int_{0}^{t}g_x(\tau)d\tau,\hspace*{1.5cm} \text{    }  \overline{G}_x(t)=\int ^{t}_{0}\overline{g}_x(\tau)d\tau \hspace*{0.5cm}
              \end{equation}
              with $\overline{g_x}(t) = \sup \left\lbrace s : g_x(s) \leqslant t\right\rbrace$, 
                  then we obtain complementary Musielak function which define their corresponding Musielak spaces $L_{G_x}$ and $L_{\overline{G}_x}$. Moreover, we suppose that\\  
              $$ (f_2): \label{f2}   \hspace{1.7cm} 1<g^-=\inf_{t\geqslant 0}\dfrac{tg_x(t)}{G_x(t)}\leqslant g^+=\sup_{t\geqslant 0}\dfrac{tg_x(t)}{G_x(t)}<+\infty. \hspace{2.3cm} \text{.}$$

                     $$ (f_3): \label{f3} \hspace{2.3cm} \lim_{t\rightarrow \infty} \dfrac{G_x(kt)}{\widehat{\varPhi}_{x,s_2}(x,t)}=0 \text{  for all  } k>0. \hspace{4.2cm} $$                                     
     $(f_4)$\label{f4} : There exist  $\theta>\varphi^+$   and  $r>0$ such that  $\text{for all } |t| \geqslant r \text{ and a.e. } x\in \Omega$         
                  $$  tf(x,t) \geqslant\theta F(x,t) \geqslant 0,            
                   $$
 where $ F(x,t):=\displaystyle\int_{0}^{t}f(x,\tau)d\tau$.\\     
          $$ (f_5):~~~~~~ \label{f5} \hspace{2.8cm} \limsup_{t\rightarrow 0} \dfrac{F(x,t)}{\varPhi_{x,y}(t)}<\dfrac{1}{\lambda_1}  \text{  uniformly for a.e. } x\in \Omega, ~~~~~~\hspace{0.7cm} $$   
where $\lambda_1$ is as in $(\ref{lambda1})$.
                 
\begin{rem}
In Problem \hyperref[P]{$(\mathcal{P}_{a})$}, since $0<s_2\leqslant s_1<1$, then by Proposition $\ref{pro21}$, $W^{s_1}_0L_{\varPhi_{x,y}}(\Omega)$ is continuously embedded in $W^{s_2}_0L_{\varPhi_{x,y}}(\Omega)$. Thus a solution for a problem of type \hyperref[P]{$(\mathcal{P}_{a})$} will be sought in $W^{s_1}_0L_{\varPhi_{x,y}}(\Omega)$.
\end{rem}      
         
   To simplify the notations, for any $i=1,2$, we set
    $$D^{s_i}u:=\dfrac{u(x)-u(y)}{|x-y|^{s_i}}, ~~~~ \Psi_{s_i}(u):= \int_{\Omega} \int_{\Omega} \varPhi_{x,y}\left( \dfrac{ |u(x)- u(y)|}{|x-y|^{s_i}}\right) \dfrac{dxdy}{|x-y|^N},  $$ 
    $$ d\mu= \dfrac{dxdy}{|x-y|^N},~~\text{ and }~~ ||u||_{s_i}:=||u||_{W^{s_i}_0L_{\varPhi_{x,y}}(\Omega)}. $$       
         Next, we give the definition of weak solutions for problem \hyperref[P]{$(\mathcal{P}_{a})$}.
   \begin{defini}
   $u\in W^{s_1}_0L_{\varPhi_{x,y}}(\Omega)$ is called a weak solution of problem \hyperref[P]{$(\mathcal{P}_{a})$} if, 
  {\small \begin{equation}\label{14}
      \int_{\Omega} \int_{\Omega} a(|D^{s_1}u|)D^{s_1}u D^{s_1}v d\mu+ \int_{\Omega} \int_{\Omega} a(|D^{s_2}u|)D^{s_2}u D^{s_2}v d\mu=\int_{\Omega}f(x,u)vdx,
         \end{equation}}
   for all $v\in W^{s_1}_0L_{\varPhi_{x,y}}(\Omega)$.
   \end{defini}
   Now, we are ready to state our existence  result.
       \begin{thm}\label{2.1.}
    Suppose that  \hyperref[f1]{$(f_1)$}-\hyperref[f3]{$(f_5)$} hold true.  
         Then problem \hyperref[P]{$(\mathcal{P}_{a})$} has a nontrivial weak solution $u\in W^{s_1}_0L_{\varPhi_{x,y}}(\Omega)$ which is a critical
         point of mountain pass type for the energy functional $J:W^{s_1}_0L_{\varPhi_{x,y}}(\Omega)\longrightarrow \R$ defined by
         \begin{equation}\label{14.}
         J(u)=\Psi_{s_1}(u)+\Psi_{s_2}(u) -\int_{\Omega}F(x,u)dx.
         \end{equation}
        \end{thm}
        
  Let us denote by $\Psi, I : W^{s_1}_0{L_\varPhi}(\Omega)\longrightarrow \R$  the functionals
       $$ \Psi(u)=\Psi_{s_1}(u)+\Psi_{s_2}(u)  \text{  }   \text{ and  } I(u)=\int_{\Omega}F(x,u)dx.$$ 
      \begin{rem}\label{rem1}
                 We note that the functional $J : W^{s_1}_0L_{\varPhi_{x,y}}(\Omega)\longrightarrow \R$ in $(\ref{14.})$ is well
                 defined. Indeed, if $u\in W^{s_1}_0L_{\varPhi_{x,y}}(\Omega)$, then by Proposition $\ref{pro21}$, we have $u\in W^{s_2}_0L_{\varPhi_{x,y}}(\Omega)$, thus, $\Psi(u)<\infty$. Moreover, by condition $\hyperref[f3]{(f_3)}$, we have that $u \in  L_{G_x}(\Omega)$
                 and thus $u \in L^1(\Omega)$. Hence, by the condition $\hyperref[f1]{(f_1)}$,
                 $$ |F(x,u)|\leqslant\int_{0}^{u}|f(x,t)|dt=c_0(|u|+G_x(|u|)).$$
                 It follows that
                 
                 $$I(u)=\int_{\Omega}|F(x,u)|dx<\infty.$$
                
                      \end{rem}          
 By a standard argument (see for instance \cite{3,sr5}) we have $J\in C^1\left( W^{s_1}_0L_{\varPhi_{x,y}}(\Omega), \R\right) $ and its G\^{a}teaux derivative is given by
 $$
 \begin{aligned}
 \left\langle J'(u),v\right\rangle =& \int_{\Omega} \int_{\Omega} a(|D^{s_1}u|)D^{s_1}u D^{s_1}v d\mu+ \int_{\Omega} \int_{\Omega} a(|D^{s_2}u|)D^{s_2}u D^{s_2}v d\mu\\
 &-\int_{\Omega}f(x,u)vdx,
  \end{aligned}$$             
	where $\langle . ,.\rangle$ denotes the usual duality between $\left(W^{s_1}_0L_{\varPhi_{x,y}}(\Omega), \|.\|_{s_1} \right) $ and its dual space $\left( \left( W^{s_1}_0L_{\varPhi_{x,y}}(\Omega)\right) ^{*},\|.\|_{s_1,*}\right).$\\            
  Next, we show an important lemma, namely that if the functional $J$ defined in $(\ref{14.})$ satisfies the conclusion of Theorem $\ref{2.2}$, then it has a critical point.
  
  \begin{lem}\label{2.4}
    Let  \hyperref[f1]{$(f_1)$}-\hyperref[f3]{$(f_3)$}  hold true. 
    Let $J$ be the functional defined in $(\ref{14.})$, and let $\left\lbrace u_n\right\rbrace$ be a sequence in $W^{s_1}_0L_{\varPhi_{x,y}}(\Omega)$ such that 
   \begin{equation}
   \label{i}  \hspace{0.5 cm} J(u_n)\longrightarrow c >0 \text{ , } \hspace{0.5 cm}  \hspace{0.5 cm} ||J'(u_n)||_{s_1,*}\longrightarrow 0.
   \end{equation} 
    Then there exists $u\in W^{s_1}_0L_{\varPhi_{x,y}}(\Omega)$ such that 
    $$J(u)=c, \hspace{0.5 cm} J'(u)=0.$$ 
    \end{lem}
       \noindent \textbf{Proof}. It follows from \hyperref[i]{$(3.4)$} that there exists $C > 0$ such that 
      $
      | J(u_n)| \leqslant C$ and $|\left\langle J'(u_n),u_n\right\rangle |\leqslant C||u_n||.$
 By assumptions  \hyperref[v1]{$(\varPhi_1)$}  and \hyperref[f1]{$(f_1)$}-\hyperref[f3]{$(f_3)$}, we have
 \begin{equation}\label{ff}
     0<t\varphi_{x,y}(t)\leqslant \varphi^+ \varPhi_{x,y}(t) \text{ for all } t> 0,
    \end{equation}
     \begin{equation}\label{ff2}
         0<tg_x(t)\leqslant g^+ G_x(t) \text{ for all } t> 0,
        \end{equation}  
    and
\begingroup\makeatletter\def\f@size{8}\check@mathfonts \begin{equation}\label{f}
 \left| \int_{\Omega\cap\left\lbrace |u_n|\leqslant r\right\rbrace }(F(x,u_n)-\theta^{-1}f(x,u_n)u_n)dx\right| \leqslant c_0\left[ (1+\theta^{-1})r+(1+\theta^{-1}g^+) G_x(r)\right]  \leqslant C.
 \end{equation}\endgroup
    Thus, by   (\ref{ff})-(\ref{f}) and Proposition \ref{1111}, we get
 
 \begin{equation} \label{28.}
   \small { \begin{aligned}
   C+C||u_n||&\geqslant J(u_n)-\dfrac{1}{\theta}\left\langle J'(u_n),u_n\right\rangle \\
   & \geqslant \Psi_{s_1}(u_n)+\Psi_{s_2}(u_n) -\dfrac{1}{\theta}\int_{\Omega} \int_{\Omega} \varphi_{x,y}\left(D^{s_1}u_n\right)D^{s_1}u_n d\mu\\
   &-\dfrac{1}{\theta}\int_{\Omega} \int_{\Omega} \varphi_{x,y}\left(D^{s_2}u_n\right)D^{s_2}u_n d\mu\\
   & -\left| \int_{\Omega\cap\left\lbrace |u_n|\leqslant r\right\rbrace }(F(x,u_n)-\theta^{-1}f(x,u_n)u_n)dx\right|\\
   &\geqslant \Psi_{s_1}(u_n)+\Psi_{s_2}(u_n) -\dfrac{\varphi^+}{\theta}\Psi_{s_1}(u_n) -\dfrac{\varphi^+}{\theta}\Psi_{s_2}(u_n) -C\\
   &\geqslant \left( 1-\dfrac{\varphi^+}{ \theta}\right)\left(\Psi_{s_1}(u_n)+\Psi_{s_2}(u_n)\right)  -C \\
   &\geqslant  \left( 1-\dfrac{\varphi^+}{ \theta}\right)\Psi_{s_1}(u_n)  -C \\
   &\geqslant \left( 1-\dfrac{\varphi^+}{ \theta}\right)||u_n||^\beta -C
      \end{aligned}}
    \end{equation}
    where $C$ denote various positive constants and $\beta= \left\lbrace \varphi^+, \varphi^-\right\rbrace $.  
 Hence, $\left\lbrace u_n\right\rbrace $ is bounded in $W^{s_1}_0L_{\varPhi_{x,y}}(\Omega)$. Since $W^{s_1}_0L_{\varPhi_{x,y}}(\Omega)$ is a reflexive space, we may assume
    that ${u_n}$ converges weakly to $u$ in $W^{s_1}_0L_\varPhi(\Omega)$. Further, since the embedding of
    $W^{s_1}_0L_{\varPhi_{x,y}}(\Omega)$ into   $L_{G_x}(\Omega)$ is compact, we obtain that $u_n\longrightarrow  u$  in $L_{G_x}(\Omega)$.\\   
     Then,  since $I\in C^1(W^{s_1}_0L_{\varPhi_{x,y}}(\Omega),\R)$, using \hyperref[f1]{$(f_1)$}-\hyperref[f3]{$(f_3)$} , we get $\lim\limits_{n\rightarrow \infty}I(u_n)=I(u)$ and $\lim\limits_{n\rightarrow \infty}I'(u_n)=I'(u)$ in $(W^{s_1}_0L_{\varPhi_{x,y}}(\Omega))^*$. Moreover, from (\ref{i}), we have $J'(u_n)\rightarrow 0$ in $(W^{s_1}_0L_{\varPhi_{x,y}}(\Omega))^*$. Hence, 
    \begin{equation}\label{37}
    \Psi'(u_n)\longrightarrow I'(u) \text{   in } (W^{s_1}_0L_{\varPhi_{x,y}}(\Omega))^*.
    \end{equation}
   Now, since $\Psi$ is convex, then we have
    $$\Psi(u_n)\leqslant \Psi(u)+\left\langle \Psi'(u_n),u_n-u\right\rangle .$$
    Therefore, using $(\ref{37})$, we may deduce that
    $$ \limsup_{n\rightarrow \infty} \Psi(u_n)\leqslant \Psi(u).$$
    It further follows from the convexity of $\Psi$ that it is weakly lower semicontinuous and hence
    $$\liminf_{n\rightarrow \infty}\Psi(u_n)\geqslant \Psi(u),$$
    which implies that
    $$\lim_{n\rightarrow \infty}\Psi(u_n)=\Psi(u).$$
    Thus
    $$\lim_{n\rightarrow \infty}J(u_n)=J(u).$$
    Consequently, by the uniqueness of the limit, we deduce that, $J'(u)=0$. The convexity of $\Psi$ implies that $\Psi'$ is
    monotone and hence
    $$\left\langle \Psi'(u_n),u_n-v\right\rangle \geqslant \left\langle \Psi'(v),u_n-v\right\rangle  \text{   for all }  v\in W^{s_1}_0L_{\varPhi_{x,y}}(\Omega). $$
    By $(\ref{37})$, we have 
    $$\left\langle I'(u)-\Psi'(v),u-v\right\rangle  \geqslant 0  \text{  for all } v\in W^{s_1}_0L_{\varPhi_{x,y}}(\Omega).$$
    Setting $v=u-th$ with $h\in W^{s_1}_0L_{\varPhi_{x,y}}(\Omega)$ and $t\in \R^+$, then, we get
    $$ \left\langle I'(u)-\Psi'(u-th),h\right\rangle  \geqslant 0 $$
    for all $h\in W^{s_1}_0L_{\varPhi_{x,y}}(\Omega)$. Letting $t\rightarrow 0$ and using the fact that $h$ is arbitrary in
    $W^{s_1}_0L_{\varPhi_{x,y}}(\Omega)$, we find that 
    $$ J'(u)=\Psi'(u)-I'(u)=0.$$
   It follows that $u$ is a critical point of $J$.\hspace*{6.2cm$\Box$ }            
  \subsection{On the geometry of the functional $J$}
  In this subsection, we will show that under the conditions which we have imposed on the functions $a_{x,y}$ and $f$, the geometric conditions \hyperref[G1]{$(G_1)$} and \hyperref[G2]{$(G_2)$} of Theorem $\ref{2.2}$ will hold. 
  \begin{lem}\label{2.3}
  Under the assumptions of Theorem $\ref{2.1.}$, the first
  geometric condition  \hyperref[G1]{$(G_1)$} of the mountain pass Theorem $\ref{2.2}$ hold
  for the functional $J$ defined in $(\ref{14.})$.
  \end{lem}
     \noindent \textbf{Proof.}
      For all $u\in W^{s_1}_0L_{\varPhi_{x,y}}(\Omega) \setminus \left\lbrace 0\right\rbrace $, the functional $I$ is satisfied:
    \begin{equation}\label{16}
       \begin{aligned}
          J(u)&=\Psi(u)-\int_{\Omega}F(x,u)dx\\
          &= \Psi(u) \left[ 1-\dfrac{\displaystyle\int_{\Omega}F(x,u)dx}{\Psi(u)}\right].
         \end{aligned}
    \end{equation}
       
     Using the condition \hyperref[f5]{$(f_5)$}, we have that there exists $\varepsilon \in (0,1)$ and $t_0>0$ such that 
     $$F(x,t)\leqslant \frac{1-\varepsilon}{\lambda_1}\varPhi_{x,y}(t) \text{ for all } |t|\leqslant t_0 \text{  and all } x\in \overline{\Omega}.$$
     We pose $\Omega_0:=\left\lbrace x\in \Omega :  |u(x)|  \geqslant t_0 \right\rbrace$, then we have 
     \begin{equation}\label{21.}
     \int_{\Omega} F(x,u(x))dx \leqslant \frac{1-\varepsilon}{\lambda_1}\int_{\Omega\setminus \Omega_0}\varPhi_{x,y}(u(x))dx +\int_{\Omega_0}F(x,u(x))dx.
     \end{equation}
  By $(\ref{11})$, we have 
\begin{equation}\label{22..}
\dfrac{(1-\varepsilon)\displaystyle\int_{\Omega\setminus \Omega_0}\varPhi(|u(x)|)dx}{\lambda_1\Psi(u)} \leqslant\dfrac{(1-\varepsilon)\displaystyle\int_{\Omega\setminus \Omega_0}\varPhi(|u(x)|)dx}{\lambda_1\varPsi_{s_1}(u)} \leqslant 1-\varepsilon.
\end{equation}
Next,  from \hyperref[f1]{$(f_1)$}, we have 
$$F(x,t)\leqslant c_0(|t|+ G_x(|t|)), \text{ for all } |t|\geqslant t_0 \text{ and  for a.e. }  x \in \Omega.$$
By proposition $\ref{pr1}$, it follows that
$$
\begin{aligned}
\int_{\Omega_0}F(x,u)dx& \leqslant c_0\left( ||u||_{L^1(\Omega)}+ \int_{\Omega}G_x(|u|)dx\right)\\
&\leqslant c_0\left( ||u||_{L^1(\Omega)}+ ||u||^{g^-}_{G_x}+||u||^{g^+}_{G_x}\right).
\end{aligned}
$$
     Since $W^{s_1}_0L_{\varPhi_{x,y}}(\Omega) \hookrightarrow L_{G_x}(\Omega)$ and $W^{s_1}_0L_{\varPhi_{x,y}}(\Omega) \hookrightarrow L^1(\Omega)$, we obtain
      \begin{equation}
      \int_{\Omega_0}F(x,u(x))dx \leqslant c_0c_1(||u||_{s_1}+||u||_{s_1}^{g^-}+||u||_{s_1}^{g^+}).
      \end{equation}
      Then,  for  $||u||_{s_1}\leqslant 1$, we find  
      \begin{equation}\label{22}
            \int_{\Omega_0}F(x,u(x))dx \leqslant 3c_0c_1||u||_{s_1},
            \end{equation} 
  where $c_1$ denote various positive constants.   By Proposition $\ref{1111}$, we have
      \begin{equation}\label{22.}
      \dfrac{\displaystyle\int_{\Omega_0}F(x,u)dx}{ \Psi_{s_1}(u)}\leqslant  3c_0c_1||u||_{s_1}^{1-\varphi^+}.
      \end{equation}
      Now, using $(\ref{16})$, $(\ref{21.})$,  $(\ref{22..})$  and $(\ref{22.})$, we obtain that 
      $$
      \begin{aligned}
      J(u)&\geqslant  \Psi(u)\left( \varepsilon -3c_0c_1||u||_{s_1}^{1-\varphi^+} \right)\\
      & \geqslant\dfrac{\varepsilon}{2} \Psi(u)
      \end{aligned}
      $$
      whenever
  \begin{equation}\label{r23}
  \rho \leqslant \min\left\lbrace  1,\left( \dfrac{\varepsilon}{6c_0c_1}\right) ^{\frac{1}{1-\varphi^+}}\right\rbrace.
  \end{equation}
  Finally, by Proposition $\ref{1111}$, we get 
  $$ ||u||\longrightarrow 0 \Longleftrightarrow \Psi(u) \longrightarrow 0.$$
  Hence for $ \rho > 0$ as given in $(\ref{r23})$, there exists a $\alpha=\alpha(\rho)>0$  such that for
  all $u$ with $||u||_{s_1}=\rho$, we have  $$\Psi(u)\geqslant \alpha.$$ We therefore
  obtain 
  $$J(u)\geqslant \alpha\dfrac{\varepsilon}{2}. $$
  Thus, if we set $r=\alpha\dfrac{\varepsilon}{2}$, we obtain that \hyperref[G1]{$(G_1)$} is satisfied.
               \hspace*{3.1cm$\Box$ }    
  \begin{lem}\label{2.3.}
    Under the assumptions of Theorem $\ref{2.1.}$, the second
    geometric condition \hyperref[G2]{$(G_2)$} of the mountain pass Theorem $\ref{2.2}$ hold
    for the functional $J$ defined in $(\ref{14.})$.
    \end{lem}
   \noindent \textbf{Proof.} 
  First, by  \hyperref[f4]{$(f_4)$}, it follows that
  \begin{equation}
  F(x,\xi)\geqslant r^{-\theta}\min\left\lbrace F(x,r), F(x,-r)\right\rbrace |\xi|^\theta
  \end{equation}
  for all $|\xi|>r$ and a.e. $x\in \Omega$. Thus by $(\ref{r1})$ and $F(x,\xi)\leqslant \max_{|\xi|\leqslant r}F(x,\xi)$ for all $|\xi|\leqslant r$, we obtain 
  \begingroup\makeatletter\def\f@size{9.5}\check@mathfonts  \begin{equation}\label{r12}
  F(x,\xi)\geqslant r^{-\theta}\min\left\lbrace F(x,r), F(x,-r)\right\rbrace |\xi|^\theta - \max_{|\xi|\leqslant r}F(x,\xi)- \min\left\lbrace F(x,r), F(x,-r)\right\rbrace 
  \end{equation}\endgroup
  for any $\xi\in \R$ and a.e. $x\in \Omega$. \\
  From Theorem \ref{TT}, we can fix $u_0\in C^\infty_0(\Omega)$ such that $||u_0||_{s_1}=1$ and let $t\geqslant 1$. By $(\ref{r12})$, we have 
{\small$$
      \begin{aligned}
      J(tu_0)&= \Psi_{s_1}(t u_0)+\Psi_{s_2}(t u_0)-\int_{\Omega}F(x,tu_0)dx\\
      &\leqslant ||tu_0||_{s_1}^{\varphi^+}+||tu_0||_{s_2}^{\varphi^+}+||tu_0||_{s_2}^{\varphi^-}-\int_{\Omega}F(x,tu_0)dx\\
      &\leqslant ||tu_0||_{s_1}^{\varphi^+}+c^{\varphi^+}||tu_0||_{s_1}^{\varphi^+}+c^{\varphi^-}||tu_0||_{s_1}^{\varphi^-}-\int_{\Omega}F(x,tu_0)dx\\
      &\leqslant (1+c^{\varphi^+})t^{\varphi^+}+c^{\varphi^+}t^{\varphi^-}-r^{-\theta}|t|^\theta\int_{\Omega}\min\left\lbrace F(x,r), F(x,-r)\right\rbrace|u_0(x)|^\theta dx\\
      &\hspace*{1cm}  + \int_{\Omega}\max_{|\xi|\leqslant r}F(x,\xi)+ \min\left\lbrace F(x,r), F(x,-r)\right\rbrace dx.
      \end{aligned}
      $$}
      From assumptions \hyperref[f1]{$(f_1)$} and \hyperref[f5]{$(f_5)$}, we get  $0<F(x,\xi)\leqslant c_0(|r|+G_x(|r|))$ for $|\xi|\leqslant r$ a.e. $x\in \Omega$. Thus, $0<\min\left\lbrace F(x,r) , F(x,-r)\right\rbrace < c_0(|r|+G_x(|r|))$, a.e.  $x\in \Omega$. Since $\theta> \varphi^+\geqslant \varphi^-$ by assumption \hyperref[f4]{$(f_4)$}, passing to the
      limit as $t\rightarrow \infty$, we obtain that $J(tu_0)\rightarrow -\infty$. Thus, the assertion \hyperref[G2]{$(G_2)$} follows by taking $e=Tu_0$ with $T$
      sufficiently large.\hspace*{8.7cm$\Box$ }          
   \subsection{Proof of Theorem \ref{2.1.}.}
  It follows from  Lemma $\ref{2.3}$ and Lemma $\ref{2.3.}$ that the hypotheses of Theorem $\ref{2.2}$ are satisfied. So Lemma $\ref{2.4}$ implies the
  existence of a nontrivial critical point of the functional $J$ which is a  weak solutions to our problem \hyperref[P]{$(\mathcal{P}_{a})$}.\\

\end{document}